%% file: main_final.tex
\documentclass[a4paper]{elsart}
\usepackage{amsmath}
\usepackage{graphicx}
\usepackage{indentfirst}
\usepackage{verbatim}
\usepackage{algorithm}
\usepackage{algpseudocode}
\usepackage{float}
\usepackage{soul}
\usepackage{subfigure}
\usepackage{multicol}
\usepackage{listings}
\usepackage{epsfig}
\usepackage{amsfonts}
\usepackage{amssymb}
\usepackage[para,online,flushleft]{threeparttable}
\usepackage{epstopdf}
\usepackage{color}
\usepackage{epstopdf}
\usepackage{hyperref}
\usepackage[left=2.54cm, right=2.54cm, top=2.54cm, bottom=2.54cm]{geometry}
\usepackage{titlesec}
\titlespacing\section{0pt}{12pt}{12pt}
\titlespacing\subsection{0pt}{12pt}{12pt}
\titlespacing\subsubsection{0pt}{12pt}{12pt}
\usepackage{bm}
\usepackage{xcolor}
\usepackage{mathtools}
\usepackage{import}
\newtheorem{theorem}{Theorem}

\newtheorem{definition}[theorem]{Definition}

\newcommand{\splitatcommas}[1]{	\begingroup
	\begingroup\lccode`~=`, \lowercase{\endgroup
		\edef~{\mathchar\the\mathcode`, \penalty0
			\noexpand\hspace{0pt plus 1em}}		}\mathcode`,="8000 #1	\endgroup
}
\usepackage{tabulary}
\newcolumntype{L}[1]{>{\raggedright\let\newline\\\arraybackslash\hspace{0pt}}m{#1}}
\newcolumntype{C}[1]{>{\centering\let\newline\\\arraybackslash\hspace{0pt}}m{#1}}
\newcolumntype{R}[1]{>{\raggedleft\let\newline\\\arraybackslash\hspace{0pt}}m{#1}}
\usepackage{multirow}

\usepackage{lineno}
\usepackage{setspace}
\begin{document}
\begin{frontmatter}
\title{A Causality-DeepONet for Causal Responses of Linear Dynamical Systems}
\author[SMU]{Lizuo Liu}
\author[Brown]{Kamaljyoti Nath}
\author[SMU]{Wei Cai}

\address[SMU]{ Department of Mathematics, Southern Methodist University, Dallas, TX 75275, USA}
\address[Brown]{ Division of Applied Mathematics, Brown University, Providence, RI 02912, USA}
\bigskip
{\bf Suggested Running Head:}
\\
A Causality-DeepONet for Causal Responses of Linear Dynamical Systems
\\
\bigskip
{\bf Corresponding Author: }
\\
Prof. Wei Cai \\
Department of Mathematics, \\
Southern Methodist University, \\
Dallas, TX 75275\\
Email: cai@smu.edu
\begin{abstract}
In this paper, we propose a DeepONet structure with causality to represent causal linear operators between Banach spaces of time-dependent signals. The theorem of universal approximations to nonlinear operators proposed in \cite{tianpingchen1995} is extended to operators with causalities, and the proposed Causality-DeepONet implements the physical causality in its framework. The proposed Causality-DeepONet considers causality (the state of the system at the current time is not affected by that of the future, but only by its current state and past history) and uses a convolution-type weight in its design. To demonstrate its effectiveness in handling the causal response of a physical system, the Causality-DeepONet is applied to learn the operator representing the response of a building due to earthquake ground accelerations. Extensive numerical tests and comparisons with some existing variants of DeepONet are carried out, and the Causality-DeepONet clearly shows its unique capability to learn the retarded dynamic responses of the seismic response operator with good accuracy.
\end{abstract}
\begin{keyword}
	Neural network, Universal approximation theory of nonlinear operator, DeepONet, Causality-DeepONet.
\end{keyword}
% \textsl{AMS Subject classifications: 65R20, 65Z05, 78M25}
\end{frontmatter}
\maketitle
%\setstretch{1.5}
%\linenumbers
\section{Introduction}
\label{Section:Introduction}

Computing operators between physical quantities defined in function spaces have many applications in forward and inverse problems in scientific and
engineering computations. For example, in wave scattering in inhomogeneous or
random media, the mapping between the media physical properties, which can be
modelled as a random field, and the wave field is a nonlinear operator, which
represents some of the most challenging computational tasks in medical imaging,
geophysical and seismic problems. A specific example comes from earthquake safety studies of
buildings and structures, the response of structures to seismic ground accelerations
gives rise to a casual operator between spaces of highly oscillatory
temporal signals.

Structural dynamic analysis has always been one of the crucial problems in the civil engineering field. Traditionally, researchers in this field analyzing structural dynamic response focus on constructing proper mathematical models like ordinary or partial differential equations and utilizing grid-based numerical methods to solve them. The finite element method \cite{zienkiewicz2005finite} is one of the popular methods considered for the solutions along with an appropriate time integration scheme like Newmark's-Beta method \cite{chopraDynamicsStructures2011,Nathan1959NewmarkBeta}. Alternatively, system identification-based methods, as an attempt to construct a surrogate model by mapping the input signals to the output responses directly, have shown their superior capability in accelerating the computations. A comprehensive review of this approach was provided in \cite{SIReview,KERSCHEN_SI}. Meanwhile, recently learning time sequential response operator between input and output signals has been studied using recurrent neural network (RNN) \cite{RNN}, long short-term memory neural network (LSTM) \cite{LSTM}, WaveNet \cite{wavenet}, the one-step ResNet approximation \cite{QIN2019620} and the multi-step recurrent ResNet approximation \cite{QIN2019620}. The RNN and its variant LSTM are ubiquitous network structures for predicting time series in financial engineering, machine translation, and sentiment analysis and so on in the natural language processing field. In particular, LSTM has been shown to have the potential to predict building responses excited by seismic ground accelerations \cite{kimrnn2020,lstm-seismic}. The one-step ResNet approximation and the multi-step recurrent ResNet approximation, provide the approximation to the integral form of the dynamical system and have been demonstrated effective equation recovery for linear and nonlinear dynamical systems \cite{Fu_2020,QIN2019620}.

Deep neural networks (DNNs), as one of the most intuitive frameworks for model reductions with its superior ability to approximate general high dimensional functions \cite{cybenko}, have been considered recently in learning mappings whose closed forms are not known. So far, DNNs have shown much promise in solving problems from scientific and engineering computing, including initial and boundary value problems of ODEs and PDEs \cite{cai19,weinan18,han18,gk20,msnn,gk19,zhangDNN}. Soon after universal approximation theorems to functions by neural networks was proposed \cite{cybenko}, Chen \& Chen \cite{tianpingchen1995} proved that there also exists a framework that could give an universal approximations to nonlinear operators between Banach spaces. Base on this theory, the DeepONet \cite{lu2021learning} was constructed for learning
operators where trunk net functions are used as a basis and the branch net functions as a mapping from the input function to a hidden manifold. The DeepONet replaced the one-hidden layer networks in the original proposal in Chen \& Chen's paper \cite{tianpingchen1995} by two deep neural networks, which has been shown to have the potential to break the curse of dimensionality from the input space. In the meantime, other approach for learning operators based on a graph kernel network \cite{li2020graphkernel} for PDEs has also been proposed. The nonlinear operator is decomposed by composing nonlinear activation functions with a class of integral operators with a trainable kernel. The Fourier neural operator has been proposed \cite{li2020fourier} by replacing the integral operator with the Fourier transform and a trainable mask in frequency domain. Both the graph kernel neural network and the Fourier neural operator shows the capability to approximate specific operators with very good accuracy and efficiency.

In this paper, we will study the DeepONet specially for operators embodying
causality from a physical system to learn mappings with causality such as those encountered in a building seismic wave response
problems. The Causality-DeepONet will be proposed to ensure the causality of
the retarded Green's function of the underlying differential equation between
the input seismic ground accelerations and the output building responses. In addition to
the causality consideration, the time homogeneity of a dynamic system will
also be used in the design of the neural network by encoding the convolutional
nature of the retarded Green's function in the choice of the network weights.
The proposed DeepONet with built-in causality allows us to learn, accurately and with minimum requirement of training data, the mapping between the
ground accelerations and the corresponding displacements of the building at the roof
level excited by the seismic ground accelerations.

It could be noted that the term causality is also used in fields such as causal inference \cite{louizos2017causal,luo2020causal} and causal interpretability of neural network \cite{CausalInterpretability}, or applying the neural network to solve problem like causal reasoning \cite{causalityReasoning}. Those work are referring to the logical cause-effect relationships between data. However, in our setting we focus on the physical concept of temporal causality, i.e., the state of the system at the current time is not affected by that of the future, but only by its current state and past history. For this reason, we name our framework Causality-DeepONet.

The rest of the paper is organized as follows. In section \ref{Section:Problem statement}, we state the problem considered in the present study. In section \ref{Section:Background}, we give a short review of the universal approximation theory of nonlinear operator by neural networks, and its recent development DeepONet. Further, we provide a review of a multi-scale neural network introduced to handle high frequency functions and the POD-DeepONet for efficient basis functions in the trunk net of DeepONet. In section \ref{Section:Methodology}, we propose two extensions of the DeepONet, one is the multi-scale DeepONet, the other is the Causality-DeepONet. In section \ref{Section:Numerical Results}, a comparison of the results with all the mentioned frameworks will be carried out. The conclusion and future works are included in the section \ref{Section:Conclusions}.

\section{Problem Statement: Calculation of building response due to seismic load}
\label{Section:Problem statement}
The problem under study is the prediction of the dynamic response of a multi-story building due to seismic loading. The equation of motion, after a finite element type discretization, for the building due to ground motions during an earthquake could be written as a dynamic systems of differential equations \cite{chopraDynamicsStructures2011},
\begin{equation}
    \bm{M}\bm{\ddot{x}} + \bm{C}\bm{\dot{x}} + \bm{K}\bm{x} = \bm{f}(t),
    \label{eq: seismic}
\end{equation}
where $\bm{M}$, $\bm{C}$ and $\bm{K}$ are the mass, damping and stiffness matrices of the system from the finite element discretization. $\bm{f}(t)$ the applied force, for our case, is due to ground motions during an earthquake and could be written as
\begin{equation}
    \bm{f}(t) = \bm{M}\bm{\iota}\ddot{u}_g,
\end{equation}
where $\ddot{u}_g$ is the ground acceleration due to the earthquake and $\bm{\iota}$ is the influence vector. The interested reader may refer \cite{chopraDynamicsStructures2011} for more details on formulation and solution methods.
\par
\begin{figure}[H]
\centering
\subfigure{\includegraphics[width=\textwidth]{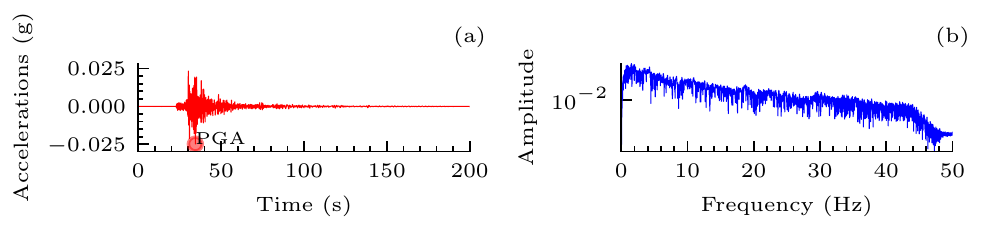}}
\caption{\textbf{A typical ground acceleration of due to earthquake:} The ground acceleration is due to 14383980 earthquake recorded at station North Hollywood, 2008 (a) time history of the acceleration (b) frequency spectrum. The absolute maximum acceleration is indicated which is also known as Peak Ground Acceleration (PGA) of the earthquake.}
\label{Figure:Typical earthquake}
\end{figure}

 Ground accelerations due to earthquakes are recorded at different recording stations. In the present study, we consider ground accelerations due to earthquakes for different earthquakes recorded at different stations and taken from the database of the Pacific Earthquake Engineering Research Center (\url{https://peer.berkeley.edu/}) \footnote{The earthquake ground acceleration considered are taken from the Pacific Earthquake Engineering Research Center (PEER: \url{https://peer.berkeley.edu/})}. One of the typical records of ground accelerations is shown in Fig. \ref{Figure:Typical earthquake}. The earthquake record at different stations may be recorded at different sampling rates (different $\delta t$). Earthquake records recorded at finer $\delta t<0.02$ sec are filtered using a butterworth filter with frequency (0.1-24.9) Hz then re-sampled to $\delta t = 0.02$ sec and after that amplified to match with original PGA level. Figures before and after processing for the mentioned earthquake record in Fig. \ref{Figure:Typical earthquake} are shown in Appendix \ref{Appendix: processed signals}. The building is considered at rest initially with the initial condition
\begin{equation}
    \left\{\begin{aligned}
       & \bm{x}(0) = \bm{0}, \\
        &    \bm{\dot{x}}(0) = \bm{0}.
    \end{aligned}\right.
    \label{eq: initial condition}
\end{equation}
Our objective is to evaluate an operator operating on the ground acceleration and predict the response of the building. Thus, it is a mapping from ground acceleration to the response of the top floor of the building.
\begin{equation}
    \mathcal{R}:\ddot{u}_g(t) \longrightarrow x_1(t).
\end{equation}
Detail numerical study is carried out for a six-storied reinforced cement concrete (RCC) building. A 3D model of the building is generated in openseespy \cite{openseespy}. Apart from the dead load (beam, column, slab, wall etc.), the live load is also considered on each floor and roof. The lumped mass of the structure is calculated for a full dead load and 50\% of the live load at floor level and 25\% at roof level. The damping matrix of the building is calculated using model damping for 5\% model damping for all the modes. Ground acceleration is applied in the major direction. The records obtained from PEER contain 3 different directions. The vertical ground accelerations are not considered. The other two horizontal ground accelerations are considered one at a time and applied only in the major direction.

\par In the case of a classical damped system \cite{chopraDynamicsStructures2011}, the displacement $\bm{x}(t)$ may be decomposed as the superposition of the modal contributions of undamped system:
\begin{equation}
\bm{x}(t)  =\sum_{l=1}^{n}\phi_{l}q_{l}\left(  t\right) \coloneqq \bm{\Phi}\bm{q},
\label{eq: displacement decomposition}
\end{equation}
where $q_l$ and $\phi_l$ are the modal coordinates and the corresponding modes, respectively, for natural frequency $\omega_l$. $\bm{\Phi}=[\phi_1,\phi_2,\cdots,\phi_n]$, $\bm{q} = [q_1, q_2,\cdots,q_n]^{T}$.
The responses (displacement) $\bm{x}_{}(t)$ of the system for ground accelerations $\ddot{u}_g(t)$ may be represented as
\begin{equation}
   \bm{x}(t) = \int_{0}^{t}  \ddot{u}_g(\tau) \bm{h}(t-\tau)d\tau,
   \label{eq: seismic-arbitrary}
\end{equation}
where
\begin{equation}
   \bm{h}(t) = -\sum_{\ell=1}^{n}\phi_\ell\dfrac{\Gamma_\ell}{\omega_{\ell D}}e^{-\xi_\ell\omega_\ell t} \sin\omega_{\ell D} t,
   \label{eq: seismic-delta}
\end{equation}
$\bm{h}(t)$ is the unit impulse response, also known as, the Green's function and fundamental solutions,  $\Gamma_\ell = \dfrac{\phi_\ell^T \bm{M} \bm{\iota}}{\phi_\ell^T\bm{M}\phi_\ell}$ and $\bm{\iota}$ is the influence vector, $\omega_{\ell D} = \omega_\ell \sqrt{1-\xi_\ell^2}$ with $\xi_\ell$ as the damping ratio.

\par In the case of a non-classical damped system \cite{chopraDynamicsStructures2011}, the responses (displacement) $\bm{x}(t)$ of the system due to ground accelerations $\ddot{u}_{g}(t)$ may be represented as
\begin{equation}
   \bm{x}_{}(t) =  -\sum_{\ell = 1}^{n} \left[ \tilde{\gamma}_{\ell}^{\delta}\omega^{\mathcal{N}}_\ell  D_{\ell}\left( t \right) + \alpha_\ell^{\delta}  \dot{D}_{\ell}\left( t \right) \right],
\label{eq: nonclassical-arbitrary}
\end{equation}
where
\begin{equation}
 D_{\ell}\left( t \right) = \int_{0}^{t}  \ddot{u}_g(\tau) H_{\ell}(t-\tau)d\tau,
\label{eq: nonclassical dn}
\end{equation}
and
\begin{equation}
    \dot{D}_{\ell}\left( t \right) = \int_{0}^{t}  \ddot{u}_g(\tau) \dot{H}_{\ell}(t-\tau)d\tau,
\label{eq: nonclassical ddn}
\end{equation}
$H_{\ell}(t) = -\dfrac{1}{\omega^{\mathcal{N}}_{\ell D}} e^{-\xi_{\ell}^{\mathcal{N}}\omega_{\ell}^{\mathcal{N}}t} \sin \omega_{\ell D}^{\mathcal{N}} t$ is the {\color{black}the Green's function of the non-classical damped system},
with
\begin{equation}
\omega^{\mathcal{N}}_{\ell D} = \omega^{\mathcal{N}}_{\ell}\sqrt{1 - (\xi^{\mathcal{N}}_{\ell})^2},
\label{eq: nonclassical omega-nd}
\end{equation}
and
\begin{equation}
\omega^{\mathcal{N}}_{\ell} = |\lambda^{\mathcal{N}}_{\ell}|,\quad \xi^{\mathcal{N}}_{\ell} = - \dfrac{\text{Re}\left( \lambda^{\mathcal{N}}_{\ell} \right)}{|\lambda^{\mathcal{N}}_{\ell}|},
\label{eq: nonclassical omega-xi}
\end{equation}
  where \(\lambda^{\mathcal{N}}_{\ell}\) is the eigenvalues  of the system of first-order differential equations reduced from Eq. $(\ref{eq: seismic})$. \(\omega^{\mathcal{N}}_{\ell}\) is a function of the amount of system damping.

 Further, \(\tilde{\gamma}_{\ell}^{\delta}=\left(\xi_{\ell}^{\mathcal{N}}\alpha_{\ell}^{\delta}-\sqrt{1 - (\xi_{\ell}^{\mathcal{N}})^2}  \gamma_\ell^{\delta} \right)\) with  \(\alpha_{\ell}^{\delta} = \text{Re}(2\beta_{\ell}^{\delta} \psi_{\ell}), \gamma_{\ell}^{\delta} = \text{Im}(2\beta_{\ell}^{\delta} \psi_{\ell})\) and
 \begin{equation}
\beta_{\ell}^{\delta} = \dfrac{-\psi_{\ell}^{T} \bm{M} \bm{\iota}}{2 \lambda_{\ell}^{\mathcal{N}}\psi_{\ell}^{T}\bm{M}\psi_{\ell}+ \psi_{\ell}^{T}\bm{C}\psi_{\ell}},
\label{eq: nonclassical bng}
\end{equation}
where $\psi_\ell$ is the corresponding eigenvector of $\lambda_\ell$ .

The discussion above of Eqs. $(\ref{eq: seismic-arbitrary})$ and \(\left( \ref{eq: nonclassical-arbitrary} \right)\) infers and inspires us to consider two phenomena in formulating an operator, the first one is that the responses at the present state is not influenced by the future ground acceleration, we understand it as causality of the system, meaning the state of system at current time should not be affected by the future, but only by its past history of the ground acceleration. The second one is the convolution nature of the Green's function kernel. As shown in many works with convolutional neural networks \cite{albawi2017cnn,o2015cnn,winovich2019convpde,zhu2019physics}, the convolution function as a specific domain knowledge for neural network to learn about the target operator. With these two insights we will construct an operator in the DeepONet framework to address both causality and convolution kernel in section \ref{Subsection:Causality DeepONet}, and name it as Causality-DeepONet.

\section{Background / Preliminary}
\label{Section:Background}
In this section, first we will discuss the universal approximation theorem for operators by neural networks. Then in section \ref{Subsection:DeepOnet} we will review the DeepONet. Multi-scale deep neural networks and POD-DeepONet will be described in sections \ref{Subsection:Multi scale} and \ref{Subsection:POD DeepOnet}, respectively.
\subsection{Universal approximation theory}
\label{Subsection:Universal approximation theory}
{\color{black}The universal approximation theory of neural networks is one of the mathematical basis for the broad applications of neural network, which justifies approximating functions or operators by weighted compositions of some functions, whose inputs are also weighted. The parameters, i.e., the weights and bias, could be obtained by minimizing a loss function with optimization algorithms such as stochastic gradient descent and its variants. Thus, neural network learning turns an approximation problem to one of optimization with respect to the parameters. In this section, we focus on the approximation of operators. }

The work of \cite{tianpingchen1995} gives a constructive procedure for
approximating nonlinear operator \(\mathcal{G}\) between continuous functions $\mathcal{G}\left( f \right)\left( x \right)$ in a compact subset of $C(\mathcal{X}) $ with \(\mathcal{X} \subseteq\mathbb{R}^{d} \) and continuous functions \(f\left( x \right)\) in a compact subset of \(C\left( \mathcal{F} \right)\) with $ \mathcal{F} \subseteq \mathbb{R}^{d}$
\begin{equation}
    \mathcal{G}: f(x)\in C(\mathcal{F}) \rightarrow \mathcal{G}(f)(x)\in
C(\mathcal{X}),
\end{equation}
where \(C\left( \mathcal{X} \right)\) and \(C\left( \mathcal{F} \right)\) are the continuous function spaces over \(\mathcal{X}\) and \(\mathcal{F}\), respectively, and \(\mathcal{X}\) and \(\mathcal{F}\) are compact subsets of \(\mathbb{R}^{d}\), the Euclidean space of dimension \(d\).
The universal approximations with respect to operators are based on the two following results:
\begin{itemize}
\item \textbf{Universal Approximation of Functions \cite{tianpingchen1995}:} Given any $\varepsilon
_{1}>0,$ there exists a positive integer \(N\), \(\left\{  \bm{w}_{k}\right\}_{k=1}^{N}  \in \mathbb{R}^{d}, \left\{ {b}_{k} \right\}_{k=1}^{N}  \in \mathbb{R}\), such that functions $\mathcal{G}(f)(x)$ selected from a compact subset \(\mathcal{U}\) of
$C(\mathcal{X})$ could be uniformly approximated by a one-hidden-layer neural network with
any Tauber-Wiener $\left(  \text{TW}\right)  $\ activation function
$\sigma_{t}$
\begin{equation}
\left|\mathcal{G}(f)(x)-\sum_{k=1}^{N}c_{k}\left(\mathcal{G}(f)\right)\sigma_{t}\left(  \bm{w}_{k}\cdot
x+b_{k}\right)  \right|\leq\varepsilon_{1},\quad\forall\text{ } x\in \mathcal{X},
\label{thm: UAT_function}
\end{equation}
where $c_{k}\left(\mathcal{G}(f)\right)$ is a linear continuous functional defined on $\mathcal{V}$
(a compact subset of $C(\mathcal{F})$), and all $\bm{w}_{k},b_{k}$ are independent of $x$ and \(f\left( x \right)\).
A Tauber-Wiener activation function is defined as follows.
\begin{definition}
Assume \(\mathbb{R}\) is the set of real numbers. $\sigma: \mathbb{R} \rightarrow \mathbb{R}$ is
called a {Tauber-Wiener $\left(  \text{TW}\right) $} function if all the linear combinations
\( g(x) = \sum\limits_{i=1}^{I}c_{i}\sigma\left(  {w}_{i}x+b_{i}\right)
\)
are dense in every $C\left[  a,b\right]  $, where \(\left\{ {w}_{i} \right\}^{I}_{i=1} , \left\{ b_{i} \right\}^{I}_{i=1} , \left\{ c_{i} \right\}^{I}_{i=1} \in\mathbb{R}\) are real constants.
\end{definition}

\medskip
\item \textbf{Universal Approximation of Functionals \cite{tianpingchen1995}:} Given any $\varepsilon_{2}>0$, there exists a positive integer \(M\), \(m\) points  \(\left\{ x_{j} \right\}_{j=1}^{m}  \in \mathcal{F}\) with  real constants \(c_{i}^{k}, {W}_{ij}^{k}, B_{i}^{k} \in \mathbb{R}, i=1,\ldots, M, j=1,\ldots ,m,\) such that a functional $c_{k}\left(\mathcal{G}(f)\right)$ could be approximated by a one-hidden-layer neural
network with any $\text{TW}$ activation
function $\sigma_{b}$
\begin{equation}
\left|c_{k}\left(\mathcal{G}(f)\right)-\sum_{i=1}^{M}c_{i}^{k}\sigma_{b}\left(  \sum_{j=1}^{m}{W}_{ij}
^{k}f\left(  x_{j}\right)  +B_{i}^{k}\right)  \right|\leq\varepsilon
_{2},\quad\forall f\in \mathcal{V}, \label{thm: UAT_functional}
\end{equation}
where the coefficients $c_{i}^{k},{W}_{ij}^{k},B_{i}^{k}$ and nodes
$\{x_{j}\}_{j=1}^{m}$ and $m,M$ are all independent of $f\left( x \right)$.
\end{itemize}

Combining these two universal approximations, the authors of \cite{tianpingchen1995} propose the universal
approximations of nonlinear operators by neural networks when restricted to the compact subset  \(\mathcal{V}\) of
the continuous function space \(C\left( \mathcal{F} \right)\) defined on a compact domain $\mathcal{F}$ in \(\mathbb{R}^{d}\). Namely, given any
$\varepsilon>0,$
 there exists positive integers \(M, N\), \(m\) points  \(\left\{ x_{j} \right\}_{j=1}^{m}  \in \mathcal{F} \subseteq \mathbb{R}^{d} \) with  real constants \(c_{i}^{k}, W_{ij}^{k}, B_{i}^{k} \in \mathbb{R}, i=1,\ldots, M, j=1,\ldots ,m,\)
 \(\left\{  \bm{w}_{k}\right\}_{k=1}^{N}  \in \mathbb{R}^{d}, \left\{ b_{k} \right\}_{k=1}^{N}  \in \mathbb{R}\) that are all independent of continuous functions $f \in \mathcal{V}\subseteq C(\mathcal{F})$ and $ x\in \mathcal{X} \subseteq \mathbb{R}^{d}$ such that
\begin{equation}
\left|\mathcal{G}(f)(x)-\sum_{k=1}^{N} \sum_{i=1}^{M}c_{i}^{k}\sigma_{b}\left(  \sum_{j=1}^{m}W_{ij}
^{k}f\left(  x_{j}\right)  +B_{i}^{k}\right)\sigma_{t}\left(  \bm{w}_{k}\cdot
x+b_{k}\right)  \right|\leq\varepsilon
\label{thm: UAT_operator}
\end{equation}

\subsection{DeepONet}
\label{Subsection:DeepOnet}
\begin{figure}[H]
\centering
\def\svgwidth{18cm}
\centering  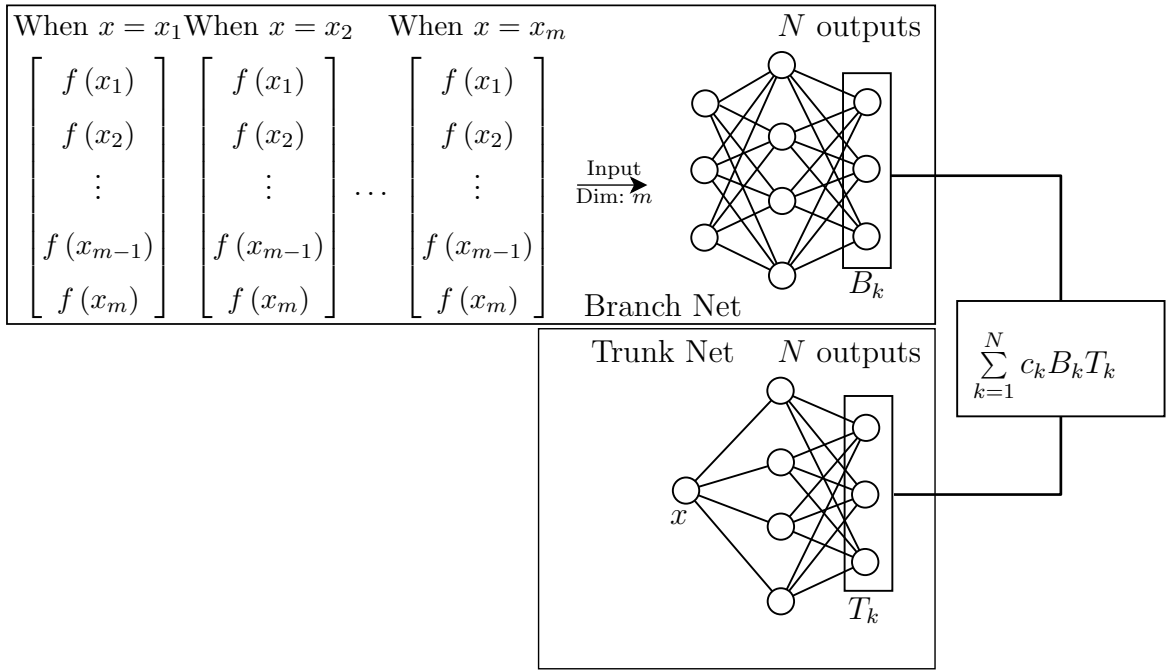
\setlength{\belowcaptionskip}{-10pt}
\caption{\textbf{Schematic Diagram of the DeepONet:} A schematic diagram of DeepONet showing branch and the trunk net along with the input data and output. The number of the input neurons of the branch net is equal to the number of sensor points in the input signals. The trunk net takes the input point $x$ at where the output function need to be evaluated. Thus, the first layer of trunk equal to the dimension of the problem. In the present study, it is the time point $x=t$ at which the output needs to be evaluated. Note for different time point \(t\), the corresponding input of branch net is the same.}
\label{fig:DeepONet}
\end{figure}

Based on the universal approximation of nonlinear operator, Lu et al. \cite{lu2021learning} proposed the deepONet by replacing the two one-hidden-layer neural networks in Eq. \(\left( \ref{thm: UAT_operator} \right)\) with two deep neural networks. For a general operator $\mathcal{G}(f)(x)$, DeepONet has form
\begin{equation}
\mathcal{G}(f)(x)\sim\sum_{k=1}^{N}c_{k}\sigma_{B,k}\left(
\left\{f\left(  x_{j}\right)\right\}_{j=1}^{m}  \right)
\sigma_{T,k}\left(   x\right),
\label{don_formula}
\end{equation}
where $\sigma_{B}\left( \cdot \right)$ with a signal $\left\{ f(x_{j}) \right\}_{j=1}^{m}$ as input is a deep neural network with \(N\) outputs, named as the branch net, $\sigma_{T}\left( \cdot \right)$ with input $x$ is also a deep neural network with \(N\) outputs is called the trunk net. The schematics are shown in Fig. \ref{fig:DeepONet}. The DeepONet has already been shown it is able to learn not only explicit mathematical operators like integration and fractional derivatives, but also PDE operators \cite{cai2021deepm,deng2022approximation,di2021deeponet,lin2021operator,lu2021learning}.

\subsection{Multi-scale  Deep Neural Network (MscaleDNN)}
\label{Subsection:Multi scale}
Multi-scale Deep Neural Network \cite{msnn} is a specific framework for problem whose output function is highly oscillatory. Since the highly oscillating feature of the responses of buildings excited by ground accelerations due to earthquakes, we also propose the multi-scale DeepONet that incorporates the multi-scale neural network into the DeepONet. The general fully connected neural network could learn the low frequency content of the data quickly, but the learning process will be stalled when higher frequency components are involved in the data. This frequency bias phenomenon is considered as the Frequency Principle studied in \cite{fprinciple}. To remedy the frequency bias in terms of learning convergence, Liu et al. \cite{msnn} introduce a MscaleDNN to accelerate the convergence of neural network  for fitting problems of highly oscillating data.
The MscaleDNN contains several sub-neural networks, for each of which a different frequency scaling is introduced by scaling the inputs accordingly, to arrive at the following form for the MscaleDNN,
\begin{equation}
f_{{\bm\theta}}\left( x \right) \sim \sum_{i=1}^{S}w_{i} f_{\bm{\theta}_{i}}\left( \mathcal{S}_{i} x \right),
\label{eq: msnn}
\end{equation}
where $\mathcal{S}_{i}$ are the custom scales and \(\left\{ f_{\bm{\theta}_{i}}\left( \cdot \right) \right\}_{i=1}^{S}\) are the distinct sub-fully connected neural networks. Eq. \(\left( \ref{eq: msnn} \right)\) shows a multi-scale neural network with \(S\) scales. The final output of multi-scale deep neural network is the weighted sum of the outputs of the sub-neural networks with trainable weights \(w_{i}\). The MscaleDNN has already shown its power for solving fitting problems and PDEs with high frequencies in \cite{linearizedLearning,msnn,msnn-stokes}.

\subsection{POD-DeepONet}
\label{Subsection:POD DeepOnet}
Instead of modeling the basis of output data by training the trunk net, Lu et al. \cite{luPOD} propose the POD-DeepONet based on the work of Bhattacharya et al. \cite{POD}. The trunk net in the vanilla DeepONet is replaced by the basis obtained from proper orthogonal decomposition (POD) of the outputs of training data after the mean of which is removed. Thus, the outputs of branch net are the coefficients of the precomputed basis vectors
\begin{equation}
\mathcal{G}(f)\sim \sum_{k=1}^{p} \sigma_{B,k}\left( \left\{f\left(  x_{j}\right)\right\}_{j=1}^{m} \right) \mathcal{B}_{k} + \mathcal{B}_{0},
\label{poddeeponet}
\end{equation}
where \(\mathcal{B}_{0}\) is the mean of the output of training data, \(\left\{ \mathcal{B}_{k} \right\}_{k=1}^{p}\) are the selected \(p\) basis vectors obtained by SVD or POD of the zero-mean outputs of training data, and \(\sigma_{B}(\left\{f\left(  x_{j}\right)\right\}_{j=1}^{m})\) is the deep neural network with  \(p\) outputs whose \(k^{th}\) output corresponds to the \(k^{th}\) singular value.
In \cite{luPOD}, it was shown that POD-DeepONet is more effective than the vanilla DeepONet and the vanilla Fourier Neural Operator \cite{li2020fourier}.

\section{Methodologies}
\label{Section:Methodology}
\subsection{Multi-scale DeepONet}
\label{Subsection:Multi scale DeepOnet}
As shown in Fig. \ref{Figure:Typical earthquake}, the spectrum of a typical earthquake signals contains not only low frequency components, but also high frequency components, therefore, we introduce the multi-scale DeepONet to handle the oscillatory information. Since the oscillatory features are function of time $t$, it is natural that we replace the fully connected trunk net by the MscaleDNN with \(S\) scales.
\begin{equation}
\begin{aligned}
    \mathcal{G}(f)(t)\sim&\sum_{k=1}^{N}c_{k}\sigma_{B,k}\left(
\left\{f\left(  t_{j}\right) \right\}_{j=1}^{m}\right)
\sigma_{T, k}\left(   t\right),\\
        \sigma_{T}\left( t \right) = & \sum_{i=1}^{S}w_{i} \sigma_{\bm{\theta}_{i}}\left( \mathcal{S}_{i} t \right),
\end{aligned}
\label{msdon_formula}
\end{equation}
where $\left\{ \sigma_{\bm{\theta}_{i}}\left( \cdot \right) \right\}_{i=1}^{S}$ are \(S\) distinct fully connected neural networks with \(N\) outputs.
\subsection{{Causality-DeepONet}}
\label{Subsection:Causality DeepONet}

The DeepONet proposed in \cite{lu2021learning} is based on a proven theorem of universal
approximation for nonlinear operators, as introduced in the section \ref{Subsection:Universal approximation theory}.
 We will apply the universal approximation theory to a nested subspaces of continuous functions indexed by the output time, which will provide a heuristic argument for the form of the Causality-DeepONet to be proposed. Rigorous mathematical justification though is still to be derived.

For a {\color{black}ground acceleration} dynamic excitation $\ddot{u}_{g}(s)$, the response function $\mathcal{R}(\ddot{u}_{g})(t)$
experiences a retardation effect due to the causality of the physical process.
Therefore, applying the universal approximation of function
(\ref{thm: UAT_function}) to input function space
\begin{equation}
C[0,t]\subseteq C[0,T],\forall t\in [0,T], \label{thm: Csubspace}
\end{equation}
we have,
\begin{equation}
    \left| \mathcal{R}(\ddot{u}_{g})(t^{\prime})-\sum_{k=1}^{N}c_{k}\left(\mathcal{R}(\ddot{u}_{g})  ,t\right)
\sigma_{t}\left(  \bm{w}_{k} t^{\prime}+b_{k}\right)
\right| \leq\varepsilon_{1},\text{ \ }\forall\text{ }t^{\prime}\in
[0,t]\subseteq[0,T].
\label{thm: UAT_function_retard}
\end{equation}

Comparing with the original universal approximation theorem of functions Eq. (\ref{thm: UAT_function}), we will require the dependence of $t$ for the parameters $c_{k}\left(\mathcal{R}(\ddot{u}_{g})  ,t\right)$, in which we may introduce the \textbf{causality} and \textbf{convolution}. It can be assumed that for every given interval \(\left[ 0,t \right] \), the approximation Eq. \(\left( \ref{thm: UAT_function_retard}\right)\) is valid within the interval based on the proof in \cite{tianpingchen1995}. In principle, the integer $N$, the parameters $\left\{ \bm{w}_{k} \right\}_{k=1}^{N}, \left\{ b_{k} \right\}_{k=1}^{N}   \in \mathbb{R},$
should also have a $t$-dependence, however, due to the nested property in
(\ref{thm: Csubspace}), for practical implementations they will be taken as global parameters to be trained for
all $t \in [0,T]$.

Following the discussion in section \ref{Section:Problem statement}, we extend the universal approximation of functionals Eq. (\ref{thm: UAT_functional}) to approximate the functionals with causality. The functional $c_{k}\left(\mathcal{R}(\ddot{u}_{g})  ,t\right)$ (from a compact subset of
$C[ 0,t] $) could be rewritten as
\begin{equation}
     c_{k}\left(\mathcal{R}(\ddot{u}_{g})  ,t\right)=c_{k}\left(  \mathcal{R}(\ddot{u}_{g} \chi_{[0,t]})
    \right)  \in\mathbb{R},\text{ \ }\ddot{u}_{g}\in C[0,t],\text{ \ }t\in\left[  0,T\right],
\label{retarded_ck}
\end{equation}
where \(\chi_{[0,t]}\left( s \right)\) is the characteristic function, such that
\begin{equation}
    \chi_{[0,t]}\left( s \right) =
    \left\{\begin{aligned}
            &1\quad s \in [0,t],\\
&0\quad \text{otherwise.}
    \end{aligned}\right.
\label{eq: characteristic}
\end{equation}
Then, following the similar approach as universal approximation of functionals Eq. \(\left( \ref{thm: UAT_functional} \right)\), we may state that
given any $\varepsilon_{2}>0$, there exists a positive integer \(M\), \(m\) {equal-spaced} points  \(\left\{ s_{j} \right\}_{j=1}^{m}  \in [ 0,t ] \) with  real constants \(c_{i}^{k}, W_{ij}^{k}, B_{i}^{k} \in \mathbb{R}, i=1,\ldots, M, j=1,\ldots ,m,\) such that \(c_{k}\left(\mathcal{R}(\ddot{u}_{g}), t \right)\) could be approximated by a one-hidden-layer neural network with any
$\text{TW}$ activation function $\sigma_{b}$
\begin{equation}
    \left\vert c_{k}(\mathcal{R}(\ddot{u}_{g}),t)-\sum_{i=1}^{M}c_{i}^{k}\sigma_{b}\left(  \sum
_{j=1}^{\left\lfloor \frac{t}{h}\right\rfloor }W_{i,m-\left\lfloor \frac{t}
{h}\right\rfloor +j}^{k}\ddot{u}_{g}\left(  s_{j}\right) +B_{i}^{k} \right)  \right\vert
\leq\varepsilon_{2},\forall \ddot{u}_{g}\in C[0,t], \label{thm: UAT_functional_causality}
\end{equation}
where the \(h\) is the step size and \(m = \left\lfloor \dfrac{t}{h} \right\rfloor \). {\color{black} For any given interval \(\left[ 0,t \right] \), based on the results of \cite{tianpingchen1995}, there exists \(m \in \mathbb{N}\) such that there are \(m\) points \(\left\{ s_{j} \right\}_{j=1}^{m}\) that could be applied to construct the functional approximation Eq. \((\ref{thm: UAT_functional})\). We further assume those points are equal-spaced. The equal-spaced sampling could be obtained by applying some appropriate smoothing kernel to input and output functions for the problems that are not equal-spaced. And likewise, the coefficients $c_{i}^{k},W_{ij}^{k},B_{i}^{k}$ and nodes
$\{s_{j}\}_{j=1}^{m}$ and $m,M$ are all independent of $\ddot{u}_{g}\left( s \right)$, however, $t$-dependent.}

The indicator function $\chi_{\lbrack0,t]}\left(  s\right)$ Eq. \(\left( \ref{eq: characteristic} \right)\), implemented
by the inner upper summation limit $\left\lfloor \frac{t}{h}\right\rfloor$ in Eq. (\ref{thm: UAT_functional_causality}), as
a discontinuous function does not belong to the continuous function space, which
could be replaced by a smoothed version with a short transition at $s=t$
while still keeping the causality.

\begin{figure}[H]
\def\svgwidth{18cm}
\centering  \import{Fig}{CausalityDeepONet.pdf_tex}
\setlength{\belowcaptionskip}{-30pt}
\caption{\textbf{Schematic Diagram of the Causality-DeepONet:} A schematic of the Causality-DeepONet showing branch and the trunk net along with the input data and output. Similar to DeepONet, the number of input neurons of branch of Causality-DeepONet is equal to the number of sensor points in the input signals. The input signals of branch, however, will be replaced by a zero-padding signals with a shifting window to express the causality and the convolution.}
\label{fig:cdeeponet}
\end{figure}
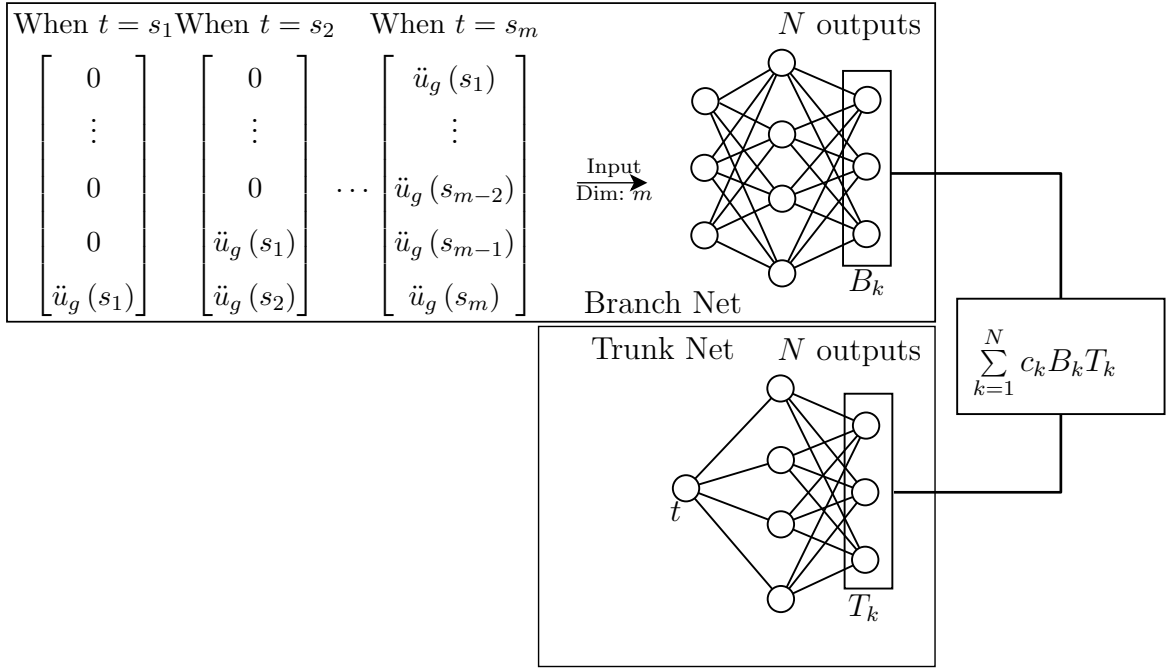

\medskip
\textbf{Causality-DeepONet:} Combining these two desired universal approximations, we can heuristically consider the following
DNN representation of an operator for retarded response for  $t\in [0,T]\subset \mathbb{R}$.
{\color{black} The basic idea is that we could find $m$ points $\{s_i\}_{i=1}^{m} \in [0,T]$ to approximate the functional $c_k(\mathcal{R}(\ddot{u}_{g}),T)$ based on the universal approximation of functionals with causality Eq. (\ref{thm: UAT_functional_causality}).  The information to approximate the functional $c_k(\mathcal{R}(\ddot{u}_{g}), t)$ where $[0,t] \subseteq [0,T]$ is offered by the value of $\{\ddot{u}_g(s_i)\}_{i=1}^{\lfloor \frac{t}{h}\rfloor }$ only. To keep the input signals of the same length at different time point, we consider zero-padding $\{\ddot{u}_g(s_i)\}_{i=1}^{\lfloor \frac{t}{h}\rfloor }$} as shown in Fig. \ref{fig:cdeeponet}. Thus the coefficients $c_{i}^{k},W_{ij}^{k},B_{i}^{k}$ and nodes $\{s_{j}\}_{j=1}^{m}$ are \(t\)-independent.
In addition, the convolution with respect to the input signals could be implemented by shifting the signals, as shown in Fig. \ref{fig:cdeeponet}.

Namely, we could find positive integers \(M, N\), \(m\) equal-spaced points  \(\left\{ t_{j} \right\}_{j=1}^{m}  \in \left[ 0,T \right] \) with  real constants \(c_{i}^{k}, W_{ij}^{k}, B_{i}^{k} \in \mathbb{R}, i=1,\ldots, M, j=1,\ldots ,m,\)
 \(\left\{  \bm{w}_{k}\right\}_{k=1}^{N}  \in \mathbb{R}, \left\{ b_{k} \right\}_{k=1}^{N}  \in \mathbb{R}\) that are all independent to continuous functions $\ddot{u}_{g} \in C\left[ 0,T \right]$ and $t$, such that
\begin{equation}
     \mathcal{R}(\ddot{u}_{g})(t) \sim \sum_{k=1}^{N}\sum_{i=1}^{M}c_{i}^{k}\sigma
    _{b}\left( \sum_{j=1}^{\left\lfloor \frac{t}{h}\right\rfloor }W_{i,\left\lfloor \frac{T}{h} \right\rfloor -\left\lfloor \frac{t}{h}\right\rfloor +j}^{k}\ddot{u}_{g}\left(  s_{j}\right)+\sum_{j=\left\lfloor \frac{t}{h}\right\rfloor }^{\left\lfloor \frac{T}{h}\right\rfloor -1 }W_{i,\left\lfloor \frac{T}{h} \right\rfloor -j}^{k} 0 + B_{i}^{k}
\right)  \sigma_{t}\left(  \bm{w}_{k} t+b_{k}\right).
\label{eq: UAT_operator_causal}
\end{equation}

\subsection{Loss function and error calculation}
\label{Subsection:Loss function and error calculation}
In the present study, we consider two loss functions depending on the method considered. Assuming $\hat{\bm{x}}_\ell(\bm{\theta})$ are the predicted response by neural network for the $\ell^\text{th}$ earthquake ground acceleration $\ddot{\bm{u}}_g^{(\ell)}$ with the corresponding true value $\bm{x}_\ell$, where \(\bm{\theta}\) are trainable variables of the neural networks, including the weights and bias. The first loss function considered is the MSE loss function
\begin{equation}
\begin{aligned}
    \mathcal{L}(\bm{\theta}) & = \dfrac{1}{n} \sum_{\ell=1}^n \left\|\bm{x}_{\ell} - \hat{\bm{x}}_\ell(\bm{\theta})\right\|^2
		\label{Eq:loss deeponet 1} \\
  & = \dfrac{1}{n}\sum_{\ell=1}^n\left[\dfrac{1}{m}\sum_{i=1}^m(x_\ell^{(i)} - \hat{x}_\ell^{(i)}(\bm{\theta}))^2 \right].
\end{aligned}
\end{equation}
where $n$ is the number of samples (different earthquake accelerations) considered and \(\|\cdot\|^2\) is the MSE error for one sample, $m$ is the number of points in each of the earthquake acceleration.

The second loss function considered is a weighted MSE loss function and defined as,
\begin{equation}
\begin{aligned}
    \mathcal{L}(\bm{\theta}) = &\frac{1}{n}\sum_{\ell=1}^{n}\frac{1}{\text{max}|\bm{x}_{\ell}|}\left\|\bm{x}_{\ell} - \hat{\bm{x}}_\ell(\bm{\theta})\right\|^2 \\
    = &\frac{1}{n}\sum_{\ell=1}^{n}\frac{1}{\text{max}|\bm{x}_{\ell}|}\left[\dfrac{1}{m}\sum_{i=1}^m(x_\ell^{(i)} - \hat{x}_\ell^{(i)}(\bm{\theta}))^2 \right].
\end{aligned}
\label{Eq:loss deeponet 2}
\end{equation}

The penalty is set to be the reciprocal of the maximum of the absolute value of the response. This act as a normalization factor when the responses have different magnitude for different earthquakes. The larger penalty is considered to the responses whose magnitude is smaller, thus it is expected that the neural network could predict the response whose magnitude is smaller accurately.

In order to check the accuracy of the predicted results we consider relative $L_2$ error,
\begin{equation}
\begin{aligned}
    \text{Relative } L_2 \text{ Error}&= \dfrac{1}{n} \sum_{\ell=1}^{n}\dfrac{\left\|\bm{x}_{\ell} - \hat{\bm{x}}_{\ell}({\bm{\theta}})\right\|}{\left\|\bm{x}_{\ell}\right\|}\\
                        &= \dfrac{1}{n} \sum_{\ell=1}^{n}\sqrt{\dfrac{\sum\limits_{i=1}^m(x_\ell^{(i)} - \hat{x}_\ell^{(i)}(\bm{\theta}))^2}{\sum\limits_{i=1}^m(x_\ell^{(i)})^2}},
\end{aligned}
    \label{Eq:RL2}
\end{equation}
and for the error with respect to the \(\ell^{th}\) case, we consider the relative error,
\begin{equation}
    \text{Err}  =  \dfrac{ \text{max}_{i}\left|x^{(i)}_{\ell} - \hat{{x}}^{(i)}_{\ell}({\bm{\theta}}) \right|}{\text{max}_{i}\left|x^{(i)}_{\ell}\right|}.
    \label{Eq:RE}
\end{equation}
The parameters $\bm{\theta}$ which include both weights and biases are optimised using the Adam optimizer \cite{Kingma_2014_adam} in a Pytorch \cite{pytorch} environment,
\begin{equation}
	\bm{\theta}^* = \text{arg\,}\underset{\bm{\theta}}{\min}\; \mathcal{L}(\bm{\theta}).
\end{equation}
Once the optimized parameters of the networks (weights and biases) are obtained, these may be used for the prediction of the response of the system for an unknown input signal (earthquake ground acceleration).

\section{Numerical Results and Discussion}
\label{Section:Numerical Results}

In this section, we will present the numerical results of multi-scale DeepONet (MS-DeepONet) and Causality-DeepONet for the prediction of response of the multistoried building discussed in section \ref{Section:Problem statement}. We will also have a comparison study of the results with a few other DeepONet methods. First, we will study the prediction of the response with different DeepONet methods along with different sizes of networks and training samples. Then, we will present predicted responses with multi-scale DeepONet and Causality-DeepONet.
We also study the different methods with different network sizes and training samples, which are discussed in subsequent sections.
To avoid overfitting, dropout \cite{Hinton2012dropout} is considered during training for few of cases, but is disabled while evaluation.

The testing dataset consists of different earthquakes which are not included in the training dataset. The test dataset is considered from 19 different earthquakes. One of them is recorded at three different stations. Two of them are recorded at the same station. Thus, the testing dataset consists of 44 ground accelerations (2 horizontal directions) from different earthquake recording stations. The training dataset consists of different earthquakes not considered in the testing dataset are may be from the same or different earthquakes and the same or different recording stations. Details about the training and testing dataset is shown in Table \ref{table:data stats} in Appendix \ref{Appendix: data statistics}.
\subsection{DeepONet and POD-DeepONet}
First, we will present the predicted results with the DeepONet method. For this purpose, we consider different trunk and branch sizes along with different training samples. As discussed in section \ref{Subsection:Loss function and error calculation}, we consider two different loss functions given by Eqs. (\ref{Eq:loss deeponet 1}) and (\ref{Eq:loss deeponet 2}).

Different network sizes considered in the branch and trunk for DeepONet are shown in Table \ref{Table:DeepONet} along with the training samples considered. The training of DeepONet is considered with Adam optimizer for a total epoch of 5000 with $\text{ReLU}(x)$ activation function with a learning rate of $10^{-4}$ in the first 1000 epochs, then $10^{-5}$ in the 1000 to 3000 epochs, and $10^{-6}$ for the remaining epochs. In order to avoid overfitting, we consider using dropout with a rate of 0.01 for the branch net and $L_2$ weight regularization with $3\times 10^{-5}$ coefficient for weights of the branch net as well. The relative $L_2$ errors for the training and testing samples after 5000 epoch are also shown in Table \ref{Table:DeepONet}. The relative $L_2$ errors with epoch for training and testing are shown in Fig. \ref{Fig:Loss DeepOnet}. A few more studies about the DeepONet are shown in the Appendix \ref{Appendix: DeepONet-further-tests}. The relative $L_2$ errors with epoch when using DeepONet with different activation functions are shown in Fig. \ref{Fig:Loss DeepOnet-extra-activation}, with $\sin(x), \tanh(x), \text{Sigmoid}(x)$ considered. The relative $L_2$ error with epoch of case that $t$ is scaled to $[0,1]$ is shown in Fig. \ref{Fig:Loss DeepOnet-extra}(a). The relative $L_2$ error with epoch of case with fixed learning rate $10^{-4}$ is shown in Fig. \ref{Fig:Loss DeepOnet-extra}(b). The relative $L_2$ errors with epoch of cases training up to 20000 epochs with fixed learning rate $10^{-4}$ is shown in Fig. \ref{Fig:Loss DeepOnet-extra}(c)-(d).  The predicted responses for few of testing dataset are shown in Appendix \ref{Appendix: predictions-DeepONet}.
It could be observed that the error in predicted responses are high for both training and testing dataset in all the cases from Table \ref{Table:DeepONet} and Fig. \ref{Fig:Loss DeepOnet} and cases in Appendix \ref{Appendix: DeepONet-further-tests}.

\begin{table}[H]
\centering
\begin{threeparttable}[b]
\setlength{\belowcaptionskip}{-10pt}
    \caption{{Relative $L_2$ error for training and testing dataset when predicted using different sizes of DeepONet.}}
    \label{Table:DeepONet}
\begin{tabular}{cccccccc}
\hline
\multirow{2}{*}{Case} & \multirow{2}{*}{Branch\tnote{1}} & \multirow{2}{*}{Trunk\tnote{1}} & \multirow{2}{*}{Sample} & \multicolumn{2}{c}{Loss Eq. (\ref{Eq:loss deeponet 1})}  & \multicolumn{2}{c}{Loss Eq. (\ref{Eq:loss deeponet 2})}\\ \cline{5-8}
& & & & Train & Test & Train & Test\\ \hline
1 & [4000]-[50]$\times$3-[50]    & [1]-[50]$\times$3-[50]   & 50 & 1.0 & 0.999 & 1.001 & 1.004\\
2 & [4000]-[100]$\times$3-[100]  & [1]-[100]$\times$3-[100] & 50 & 1.0 & 0.999 & 1.00 & 1.002\\
3 & [4000]-[200]$\times$3-[200]  & [1]-[200]$\times$3-[200] & 50 & 1.0 & 0.999 & 1.001 & 1.003\\ \hline
4 &[4000]-[50]$\times$3-[50]    & [1]-[50]$\times$3-[50]    & 100 & 1.00 & 0.999 & 1.00 &1.002\\
5 &[4000]-[100]$\times$3-[100]  & [1]-[100]$\times$3-[100]  & 100 & 1.0 & 0.999 & 1.0  &1.002\\
6 &[4000]-[200]$\times$3-[200]  & [1]-[200]$\times$3-[200]  & 100 & 1.0  & 0.999 & 1.0 &0.999\\ \hline
\end{tabular}
\begin{tablenotes}
\item [1] The notation $[N_1]$-$[N_2]\times3$-$[N_3]$ represents a neural network with the input size of $N_1$, 3 hidden layers with $N_2$ neurons in each layer, and the output dimension of $N_3$ neurons.
\end{tablenotes}
\end{threeparttable}
\end{table}
\begin{figure}[H]
\centering
\subfigure{\includegraphics[width=0.32\textwidth]{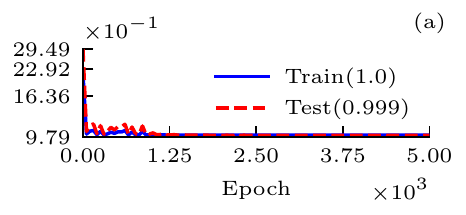}}
\subfigure{\includegraphics[width=0.32\textwidth]{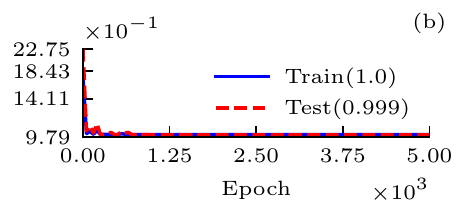}}
\subfigure{\includegraphics[width=0.32\textwidth]{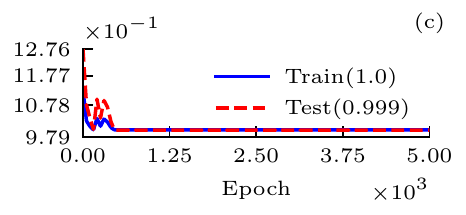}}
\subfigure{\includegraphics[width=0.32\textwidth]{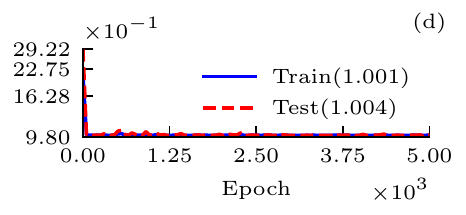}}
\subfigure{\includegraphics[width=0.32\textwidth]{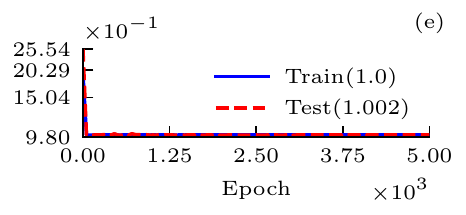}}
\subfigure{\includegraphics[width=0.32\textwidth]{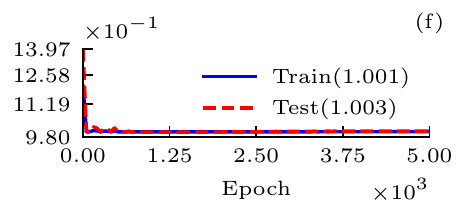}}
\subfigure{\includegraphics[width=0.32\textwidth]{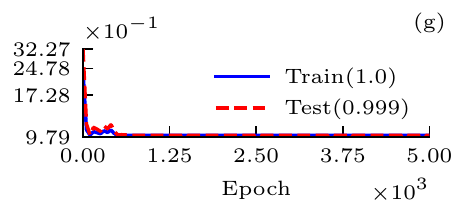}}
\subfigure{\includegraphics[width=0.32\textwidth]{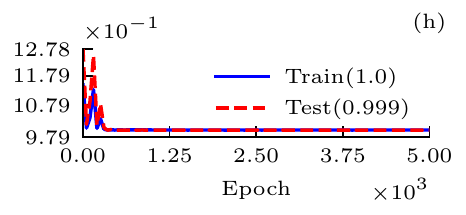}}
\subfigure{\includegraphics[width=0.32\textwidth]{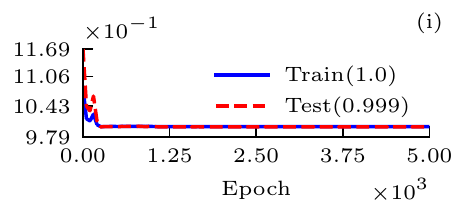}}
\subfigure{\includegraphics[width=0.32\textwidth]{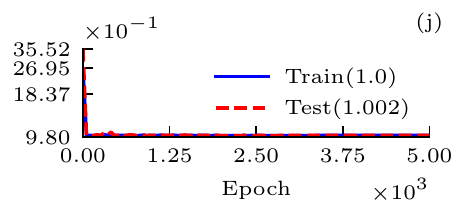}}
\subfigure{\includegraphics[width=0.32\textwidth]{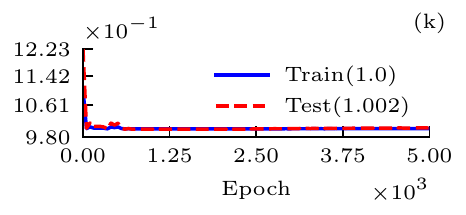}}
\subfigure{\includegraphics[width=0.32\textwidth]{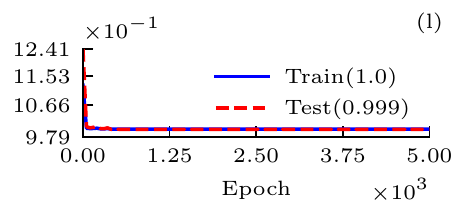}}
\setlength{\belowcaptionskip}{-10pt}
\caption{\textbf{Relative $\bm{L_2}$ Error for DeepONet:} Training and testing Relative $L_2$ error with epoch for DeepONet when considered different sample size in training and with different trunk and branch sizes. The plots (a), (c), (e), (g) (i) (k) are the cases corresponding to Case-1 to Case-6 of Table \ref{Table:DeepONet}, respectively, with Loss function Eq. (\ref{Eq:loss deeponet 1}). Similarly, the plots (b), (d), (f), (h), (j) are the cases corresponding to Case-1 to Case-6 of Table \ref{Table:DeepONet}, respectively, with Loss function Eq. (\ref{Eq:loss deeponet 2}).
}\label{Fig:Loss DeepOnet}
\end{figure}

To understand the effect of normalization of the input and output dataset on the accuracy of the predicted results, we consider Gaussian normalization of the input and output data,
\begin{equation}
	x_{\text{norm}}(t) = \dfrac{x(t) - \mu_x(t)}{\sigma_x(t)}
\end{equation}
where $x_{\text{norm}}(t)$ are the data after normalization. $\mu_x(t)$ and $\sigma_x(t)$ are the ensemble mean and standard deviation of the training dataset.

The predicted responses are decoded to the actual response with the same mean and standard deviation. Similar to the DeepONet case discussed above, we study the effect of normalization with different network sizes and the results are shown in Table \ref{Table:DeepONet normalized} along with training and testing loss for only the MSE Loss function Eq. (\ref{Eq:loss deeponet 1}) considered. A $\text{ReLU}(x)$ activation function is considered with learning rate of $10^{-3}$ in the first 1000 epochs, then $10^{-4}$ in the 1000 to 10000 epochs, $10^{-5}$ for rest of the epoch up to 20000 epochs. The other hyper-parameters considered are a dropout rate of 0.01 and a $L_2$ weight regularization with a coefficient of $10^{-5}$ for the branch net.
The relative $L_2$ error with epoch for training and testing dataset are shown in Fig. \ref{Fig:Loss DeepONet normalized}. It could be observed that the training relative $L_2$ error does not reduce even after the normalization of the input and output.

\begin{table}[H]
\centering
\begin{threeparttable}[b]
\setlength{\belowcaptionskip}{-10pt}
\caption{Relative $L_2$ error for training and testing dataset when predicted using different sizes of DeepONet with Gaussian normalization for input and output.}
\label{Table:DeepONet normalized}
\begin{tabular}{cccccc}
\hline
\multirow{2}{*}{Case} & \multirow{2}{*}{Branch} & \multirow{2}{*}{Trunk}&  \multirow{2}{*}{Sample} & \multicolumn{2}{c}{Loss Eq. (\ref{Eq:loss deeponet 1})} \\ \cline{5-6}
& & & & Train & Test \\ \hline
1 &[4000]-[100]$\times$3-[100]  & [1]-[100]$\times$3-[100]  & 100 & 1.847 & 2.379 \\
2 &[4000]-[200]$\times$3-[200]  & [1]-[200]$\times$3-[200]  & 100 & 1.841 & 2.378 \\
\hline
\end{tabular}
\end{threeparttable}
\end{table}

\begin{figure}[H]
\centering
\subfigure{\includegraphics[width=0.48\textwidth]{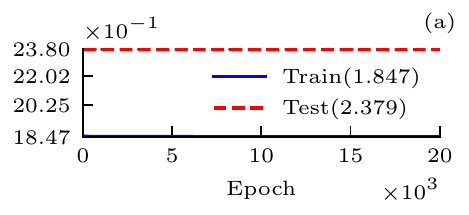}}
\subfigure{\includegraphics[width=0.48\textwidth]{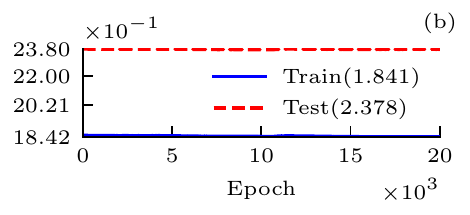}}
\setlength{\belowcaptionskip}{-10pt}
\caption{\textbf{Relative $\bm{L_2}$ Error for DeepONet with Normalization:} Relative $L_2$ Error for training and testing dataset when using different sizes of DeepONet with Gaussain normalization for input and output. (a) and (b) are Case-1 and Case-2 of Table \ref{Table:DeepONet normalized}, respectively.
}\label{Fig:Loss DeepONet normalized}
\end{figure}

From the above discussion, it could be observed that the DeepONet is not able to predict the response of the building with sufficient accuracy. As discussed in section \ref{Section:Background}, one of the modified versions of DeepONet is the POD-DeepONet, where the trunk net is replaced by POD modes obtained by proper orthogonal decomposition of the zero-mean training data (output data). These basis act as the trunk, and the branch is expected to learn the coefficient of the basis vectors.

\begin{figure}[H]
\centering
\subfigure{\includegraphics[width=0.32\textwidth]{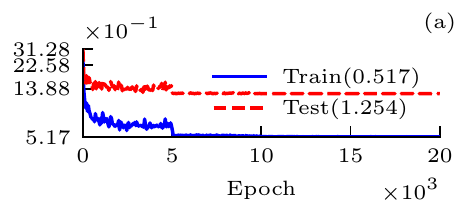}}
\subfigure{\includegraphics[width=0.32\textwidth]{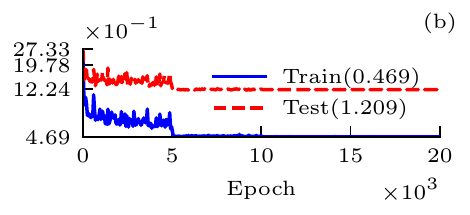}}
\subfigure{\includegraphics[width=0.32\textwidth]{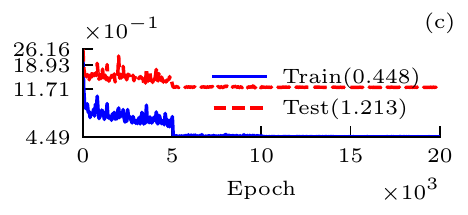}}
\subfigure{\includegraphics[width=0.32\textwidth]{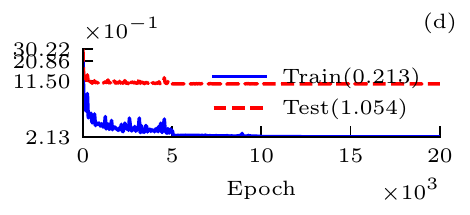}}
\subfigure{\includegraphics[width=0.32\textwidth]{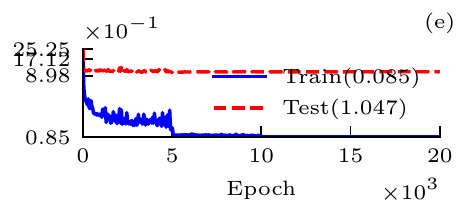}}
\subfigure{\includegraphics[width=0.32\textwidth]{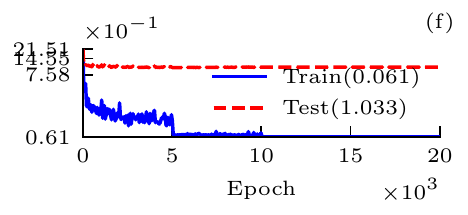}}
\caption{\textbf{Relative $\bm{L_2}$ error for POD-DeepONet:} Relative $L_2$ Error for training and testing dataset when using POD-DeepONet with different sizes. The plots (a), (c), (e) are the cases corresponding to Case-1 to Case-3 of Table \ref{Table:PODDeepONet}, respectively, with Loss function Eq. (\ref{Eq:loss deeponet 1}). Similarly, the plots (b), (d), (f) are the cases corresponding to Case-1 to Case-3 of Table \ref{Table:PODDeepONet}, respectively, with Loss function Eq. (\ref{Eq:loss deeponet 2}).
}\label{Fig:Loss PODDeepOnet}
\end{figure}

Similar to DeepONet, in this case as well, we consider two loss functions given by Eqs. (\ref{Eq:loss deeponet 1}) and (\ref{Eq:loss deeponet 2}) with different sizes of branch net as shown in Table \ref{Table:PODDeepONet}. The training of POD-DeepONet is considered with Adam optimizer for a total of 20000 epochs with $\text{ReLU}(x)$ activation function. To avoid overfitting, we consider $L_2$ weight regularization with coefficient $10^{-6}$. The learning rate considered is $10^{-3}$ in the first 5000 epochs, then $10^{-4}$ in the 5000 to 10000 epochs, and $10^{-5}$ for the rest of the epochs up to 20000. The relative $L_2$ error with epoch is shown in Fig. \ref{Fig:Loss PODDeepOnet}. It could be observed that the loss function with an additional penalty given by Eq. (\ref{Eq:loss deeponet 2}) could offer better convergence for training cases compared with the case only considering MSE loss given by Eq. (\ref{Eq:loss deeponet 1}). Further, the relative $L_2$ errors are smaller compared to DeepONet for the training dataset. However, the performance of the network is poor in the case of the testing dataset shows that 100 training samples could not offer enough information for the target response space, even though there are slight improvements with the increase in network size. The relative $L_2$ errors for training and testing dataset with epoch when using POD-DeepONet with different activation functions are shown in Fig. \ref{Fig:Loss POD-DeepOnet-extra}. The predicted response for the best training and testing cases are shown in Fig. \ref{Fig:POD-DeepOnet Pred-Best}. The performance of POD-DeepONet highly relies on the quality of the training data. If the training dataset covers a large enough region of the space of interest,
then the POD-DeepONet could have very good performance.
\begin{table}[H]
\centering
\begin{threeparttable}[b]
\caption{Relative $L_2$ error for training and testing dataset when using different sizes of POD-DeepONet.}
\label{Table:PODDeepONet}
\begin{tabular}{ccccccc}
\hline
\multirow{2}{*}{Case} & \multirow{2}{*}{Network Size\tnote{1}} &  \multirow{2}{*}{Sample} & \multicolumn{2}{c}{Loss Eq. (\ref{Eq:loss deeponet 1})}  & \multicolumn{2}{c}{Loss Eq. (\ref{Eq:loss deeponet 2})}\\ \cline{4-7}
& &  & Train & Test & Train & Test\\ \hline
1 & [4000]-[50]$\times$3-[100]    & 100 & 0.517 & 1.254 & 0.213 & 1.054\\
2 & [4000]-[100]$\times$3-[100]   & 100 & 0.469 & 1.209 & 0.085 & 1.047\\
3 & [4000]-[200]$\times$3-[100]   & 100 & 0.448 & 1.213 & 0.061 & 1.033\\ \hline
\end{tabular}
\begin{tablenotes}
\item [1]
The notation $[N_1]$-$[N_2]\times$3-$[N_3]$ represents a neural network with the input size of $N_1$, 3 hidden layers with $N_2$ neurons in each layer, and the output dimension of $N_3$ neurons.
\end{tablenotes}
\end{threeparttable}
\end{table}
\vspace*{-0.5cm}

\subsection{Multi-scale DeepONet}
In the previous section, we discussed the results of DeepONet and POD-DeepONet and observed that the results were not satisfactory. In this section, we will present and discuss the results of one of the proposed variants of DeepONet, the multi-scale DeepONet.
\begin{table}[H]
\centering
\begin{threeparttable}[b]
\setlength{\belowcaptionskip}{-10pt}	
	\caption{  Relative $L_2$ Error for Training and Testing Dataset when Using Different Sizes of Multi-scale DeepONet.}
	\label{Table:MSDeepONet}
	\begin{tabular}{cccccccc}
		\hline
    \multirow{2}{*}{Case} & \multirow{2}{*}{Branch\tnote{1}} & \multirow{2}{*}{Trunk\tnote{2}}&  \multirow{2}{*}{Sample} & \multicolumn{2}{c}{Loss Eq. (\ref{Eq:loss deeponet 1})}& \multicolumn{2}{c}{Loss Eq. (\ref{Eq:loss deeponet 2})} \\ \cline{5-8}
		& & & & Train & Test & Train & Test\\ \hline
1 &[4000]-[200]$\times$3-[200]  & [1]-20$\times$\{[10]$\times$3\}-[200]  & 100 & 0.349 & 1.041 & 0.304 & 0.994 \\
2 &[4000]-[400]$\times$3-[400]  & [1]-20$\times$\{[20]$\times$3\}-[400]  & 100 & 0.232 & 1.005 & 0.166 & 1.006 \\
\hline
	\end{tabular}
\begin{tablenotes}
\item [1]
The notation [$N_1$]-[$N_2$]$\times$3-[$N_3$] represents a neural network with the input size of $N_1$, 3 hidden layers with $N_2$ neurons in each layer, and the output dimension of $N_3$ neurons.
\item [2] The notation [$N_1$]-20$\times$\{[$N_2$]$\times$3\}-[$N_3$] represents a neural network with the input size of $N_1$, and $20$ sub-neural networks that contain 3 hidden layers with $N_2$ hidden neurons in each layer. The output dimension is $N_3$. To keep the number of neurons as the same as previous cases, the number of hidden neurons for each subnet at each layer are divided by \(20\), the number of scales.
\end{tablenotes}
\end{threeparttable}
\end{table}
\begin{figure}[H]
\centering
\subfigure{\includegraphics[width=0.45\textwidth]{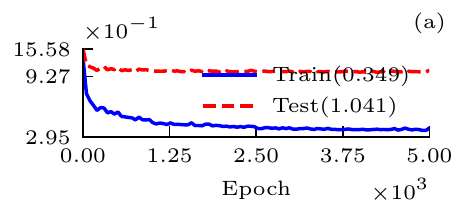}}
\subfigure{\includegraphics[width=0.45\textwidth]{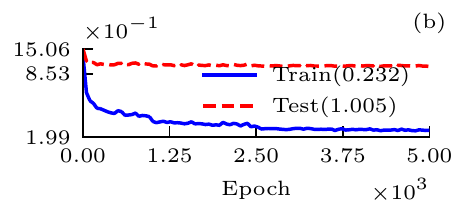}}
\subfigure{\includegraphics[width=0.45\textwidth]{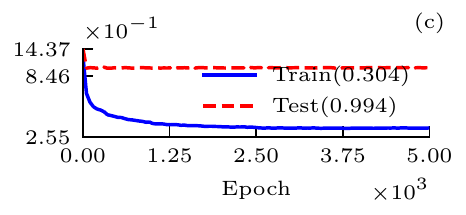}}
\subfigure{\includegraphics[width=0.45\textwidth]{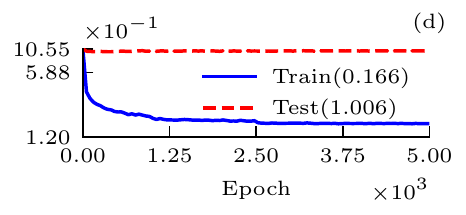}}
\setlength{\belowcaptionskip}{-10pt}
\caption{\textbf{Relative $\bm{L_2}$ Error for Multi-scale DeepONet:} Relative $L_2$ Error with epoch for  training and testing dataset when using different sizes of multi-scale DeepONet. The plots (a) and (b) correspond to Case-1 and Case-2 of Table \ref{Table:MSDeepONet}, respectively, with loss function given by Eq. (\ref{Eq:loss deeponet 1}). The plots (c) and (d) correspond to Case-1 and Case-2 of Table \ref{Table:MSDeepONet}, respectively, with loss function Eq. (\ref{Eq:loss deeponet 2}).
}\label{Fig:Loss MSDeepONet}
\end{figure}

The architecture of multi-scale DeepONet is discussed in section \ref{Subsection:Multi scale DeepOnet} where the fully connected deep neural network in the trunk net is replaced by a multi-scale neural network (MscaleDNN). In the present study we consider an MscaleDNN in the trunk with 20 equally spaced scales $[1, 1+20\pi, \cdots, 1+ 20n\pi,\cdots, 1+780\pi]$. In the meantime, the time $t$ is scaled to $t \in [0,1]$.
The activation function considered for all the cases is $\sin(x)$, according to the results in \cite{linearizedLearning,msnn-stokes}. To avoid overfitting, we consider dropout rate of 0.10 for the trunk net. The learning rate considered is $3\times10^{-4}$ in the first 1000 epochs, then $1.5\times10^{-4}$ in the 1000 to 2500 epochs, and $7.5\times10^{-5}$ for the rest of the epochs up to 5000.
The relative $L_2$ errors for different network sizes with different training loss functions are shown in Table \ref{Table:MSDeepONet}. The relative $L_2$ errors with epoch for training and testing dataset are shown in Fig. \ref{Fig:Loss MSDeepONet}. It could be observed the obtained operator is not desired based on the results of testing cases, though the MS-DeepONet accelerated the convergence for the training process.

\subsection{Causality-DeepONet}
As discussed in section \ref{Section:Problem statement} and section \ref{Subsection:Causality DeepONet}, both convolution and causality are considered in the formulation of Causality-DeepONet. In this section, we will present the numerical results and a comprehensive discussion about Causality-DeepONet for the prediction of the responses of the problem discussed in section \ref{Section:Problem statement}.

Similar to the previous numerical examples, in this study as well, we consider the two loss functions given by Eqs. (\ref{Eq:loss deeponet 1}) and (\ref{Eq:loss deeponet 2}) with different network sizes. We also study the effect of the number of training samples on the accuracy of test results. Further, given the fact that there are multiple choices of activation functions, we test a few of the popular activation functions with the same training dataset and the same network sizes. As shown later, Causality-DeepONet with standard sigmoid activation functions converges much slower. By defining a custom sigmoid function, the performance is improved. Furthermore, to study the effect of only causality without convolution, we do a numerical study with causality only and observed that convolution is also an indispensable component of the proposed Causality-DeepONet. Unlike the previous studies, we provide the initial conditions as additional data pairs \(\{\ddot{u}_g: \left[ 0,0,\ldots ,0 \right], x(0):0\} \) in the training dataset to force the Causality-DeepONet satisfy the initial conditions.

\begin{figure}[H]
\centering
\subfigure{\includegraphics[width=\textwidth]{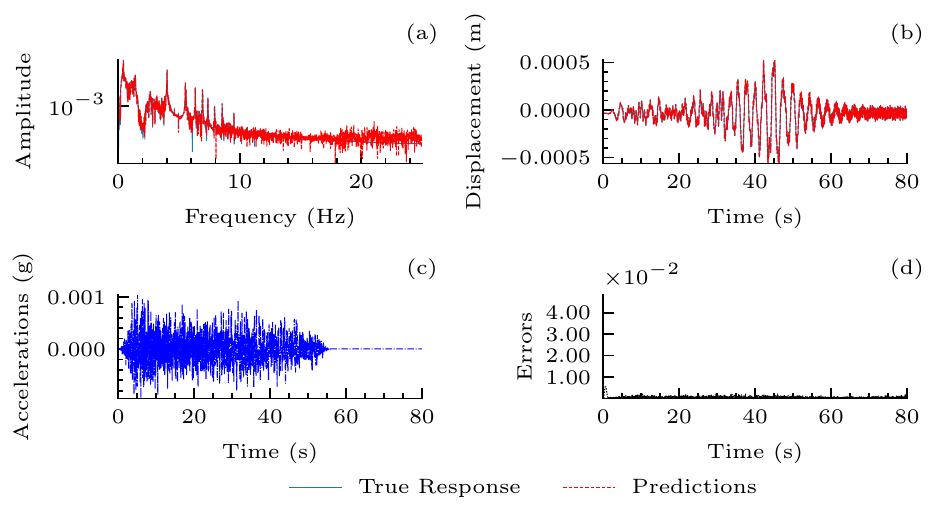}}
\setlength{\belowcaptionskip}{-10pt}
\caption{\textbf{The Worst Case of Predictions of Causality-DeepONet(Relative $\bm{L_2}$ Error: 0.0042):} The worst predictions in testing dataset for the Causality-DeepONet.
(a) The Amplitude of the prediction and true response in Fourier Domain,
(b) The prediction and true response,
(c) The corresponding input signals,
(d) The relative error Eq. (\ref{Eq:RE}).
}\label{Fig:Causality DeepOnet with bias-Pred-Worst}
\end{figure}

\begin{figure}[H]
\centering
\subfigure{\includegraphics[width=\textwidth]{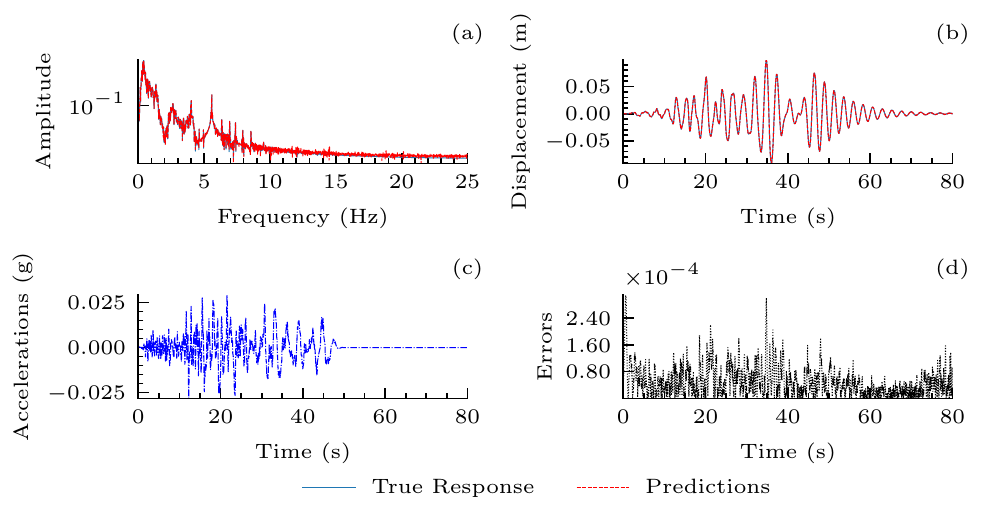}}
\caption{\textbf{The Best Case of Predictions of Causality-DeepONet(Relative $\bm{L_2}$ Error: 0.00025)} The best predictions in testing dataset for the Causality-DeepONet.
(a) The Amplitude of the prediction and true response in Fourier Domain,
(b) The prediction and true response,
(c) The corresponding input signals,
(d) The relative error Eq. (\ref{Eq:RE}).
}
\label{Fig:Causality DeepOnet with bias-Pred-Best}
\end{figure}

Fig. \ref{Fig:Causality DeepOnet with bias-Pred-Worst}  and Fig. \ref{Fig:Causality DeepOnet with bias-Pred-Best} show the worst and the best predictions using Causality-DeepONet for the test dataset. The network is trained using 100 training samples. The network size considered is [4000]-[120]$\times$2-[120] for branch and [1]-[120]$\times$2-[120] for trunk. The activation function considered is $\tanh(x)$. To avoid overfitting, we consider $L_2$ weight regularization with a coefficient $1\times 10^{-4}$ for the branch net. The learning rate considered is $10^{-3}$ in the first 2000 epochs, then $10^{-4}$ in the 2000 to 10000 epochs, and $10^{-5}$ for the rest of the training up to 20000 epochs. The loss function considered for this case is Loss Eq. (\ref{Eq:loss deeponet 2}). It could be observed that Causality-DeepONet can predict the responses with good accuracy for all the cases, as the error for the worst case also is within the satisfactory limit.


\begin{figure}[H]
\centering
\subfigure{\includegraphics[width=0.32\textwidth]{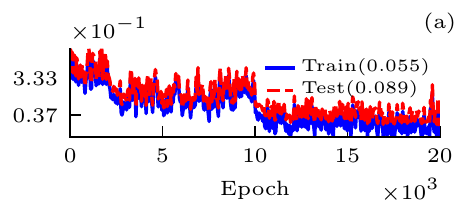}}
\subfigure{\includegraphics[width=0.32\textwidth]{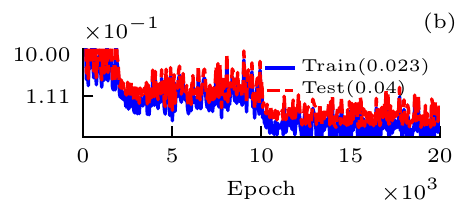}}
\subfigure{\includegraphics[width=0.32\textwidth]{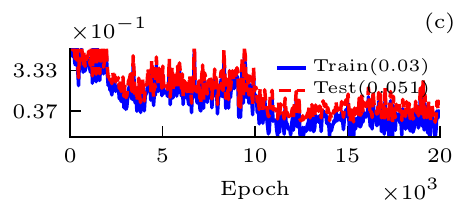}}
\subfigure{\includegraphics[width=0.32\textwidth]{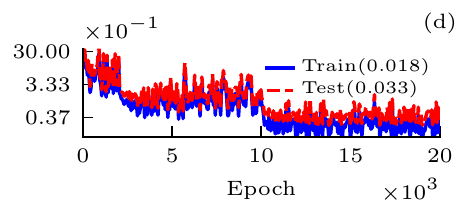}}
\subfigure{\includegraphics[width=0.32\textwidth]{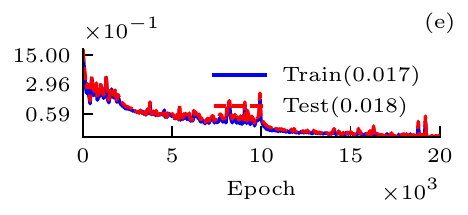}}
\subfigure{\includegraphics[width=0.32\textwidth]{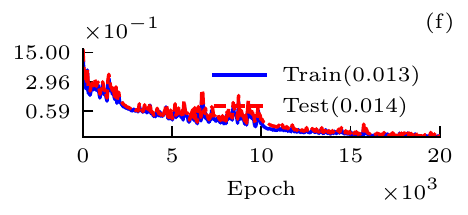}}
\subfigure{\includegraphics[width=0.32\textwidth]{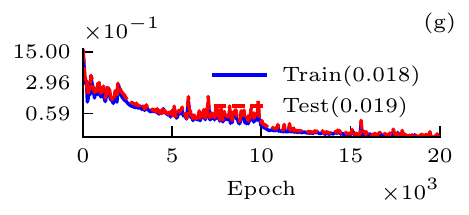}}
\subfigure{\includegraphics[width=0.32\textwidth]{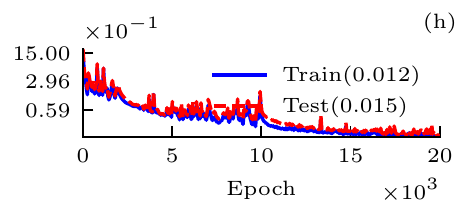}}
\caption{\textbf{Relative $\bm{L_2}$ Error for Training and Testing Dataset when Using Causality-DeepONet with Different Branch and Trunk Sizes:}  The plots (a), (b), (c), (d) are the cases corresponding to Case-1 to Case-4 of Table \ref{Table:Causality DeepONet}, respectively, with Loss function Eq. (\ref{Eq:loss deeponet 1}). Similarly, the plots (e), (f), (g), (h) are the cases corresponding to Case-1 to Case-4 of Table \ref{Table:Causality DeepONet}, respectively, with Loss function Eq. (\ref{Eq:loss deeponet 2}).
}
\label{Fig:Loss Causality DeepOnet}
\end{figure}

To study the effect of different network size and loss function (Eqs. (\ref{Eq:loss deeponet 1}) and (\ref{Eq:loss deeponet 2})), we consider different network sizes with the same number of training dataset. The activation function considered for all the cases is $\tanh(x)$. To avoid overfitting, we consider $L_2$ weight regularization for the branch net with a coefficient $10^{-4}$. The learning rate considered is $10^{-3}$ in the first 2000 epochs, then $10^{-4}$ in the 2000 to 10000 epochs, and $10^{-5}$ for the rest of the epochs up to 20000 epochs. The relative $L_2$ errors for both training and testing dataset after completion of training is shown in Table \ref{Table:Causality DeepONet} and the corresponding relative $L_2$ error with epoch is shown in Fig. \ref{Fig:Loss Causality DeepOnet}. It could be observed that the proposed Causality-DeepONet shows a good accuracy in both the case of loss functions. The MSE loss function given by Eq. (\ref{Eq:loss deeponet 1}) is more sensitive to the network size as relative $L_2$ errors for both training and testing are reduced with an increase in network sizes. The weighted loss function given by Eq. (\ref{Eq:loss deeponet 2}) is less sensitive to the network sizes for this numerical study. Further, it is also observed that the relative $L_2$ error in the case of loss function Eq. (\ref{Eq:loss deeponet 2}) is less than that of relative $L_2$ error in the case of loss function Eq. (\ref{Eq:loss deeponet 1}). Thus, we conclude that the additional penalty terms in the loss function remove the bias from the magnitude of output functions/data in the present study. For the further numerical studies conducted, we consider with loss function given by Eq. (\ref{Eq:loss deeponet 2}) only.

\begin{figure}[H]
\centering
\subfigure{\includegraphics[width=0.32\textwidth]{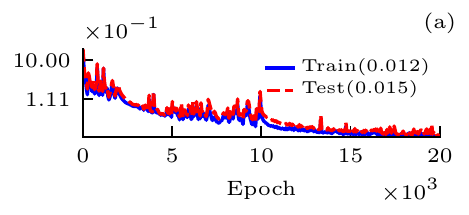}}
\subfigure{\includegraphics[width=0.32\textwidth]{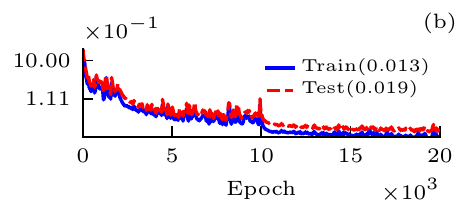}}
\subfigure{\includegraphics[width=0.32\textwidth]{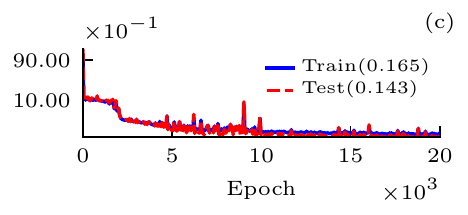}}
\subfigure{\includegraphics[width=0.32\textwidth]{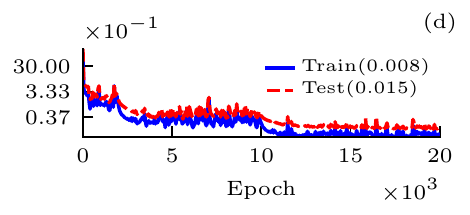}}
\subfigure{\includegraphics[width=0.32\textwidth]{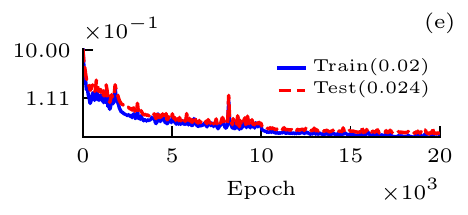}}
\caption{\textbf{Relative $\bm{L_2}$ Error for Training and Testing Dataset when Using Causality-DeepONet With Different Activation Functions:} The plots (a)-(e) correspond to cases Case-1 to Case-5 of Table \ref{Table:Causality DeepONet-activation},  respectively.
}\label{Fig:Loss Causality DeepOnet-activation}
\end{figure}

\begin{table}[H]
\centering
\begin{threeparttable}[b]
\caption{Relative $L_2$ Error for training and testing dataset when using different sizes of Causality-DeepONet}
\label{Table:Causality DeepONet}
\begin{tabular}{cccccccc}
\hline
\multirow{2}{*}{Case} & \multirow{2}{*}{Branch\tnote{1}} & \multirow{2}{*}{Trunk\tnote{1}} & \multirow{2}{*}{Sample} & \multicolumn{2}{c}{Loss Eq. (\ref{Eq:loss deeponet 1})}  & \multicolumn{2}{c}{Loss Eq. (\ref{Eq:loss deeponet 2})}\\ \cline{5-8}
& & & & Train & Test & Train & Test\\ \hline
1 & [4000]-[30]$\times$2-[30]    & [1]-[30]$\times$2-[30]   & 10 & 0.055 & 0.089 & 0.017 & 0.018\\
2 & [4000]-[60]$\times$2-[60]  & [1]-[60]$\times$2-[60] & 10 & 0.023 & 0.040 & 0.013 & 0.014\\
3 & [4000]-[90]$\times$2-[90]  & [1]-[90]$\times$2-[90] & 10 & 0.030 & 0.051 & 0.018 & 0.019\\
4 &[4000]-[120]$\times$2-[120]    & [1]-[120]$\times$2-[120]    & 10 & 0.018 & 0.033 & 0.012 & 0.015\\
\hline
\end{tabular}
\begin{tablenotes}
\item [1]
The notation [$N_1$]-[$N_2$]$\times$2-[$N_3$] represents a neural network with the input size of $N_1$, 2 hidden layers with $N_2$ neurons in each layer, and the output dimension of $N_3$ neurons.
\end{tablenotes}
\end{threeparttable}
\end{table}
\newpage
As discussed earlier, there are multiple choices of activation functions, and we test the performance of Causality-DeepONet with a few popular activation functions. For this purpose, we consider the same network sizes and training dataset and loss functions (Eq. (\ref{Eq:loss deeponet 2}) for all the cases of activation function considered. The network sizes considered are [4000]-[120]$\times$2-[120] for branch and  [1]-[120]$\times$2-[120] for trunk. To avoid overfitting, we consider $L_2$ weight regularization for the parameters in branch net with a coefficient $10^{-4}$ for case 1,2,4,5 and $5\times10^{-6}$ for case-3 in Table \ref{Table:Causality DeepONet-activation}. The learning rate considered is $10^{-3}$ in the first 2000 epochs, then $10^{-4}$ in the 2000 to 10000 epochs, and $10^{-5}$ for the rest of the epoch up to 20000 epochs.
As shown in Table \ref{Table:Causality DeepONet-activation} and Fig. \ref{Fig:Loss Causality DeepOnet-activation}, It could be concluded that the Causality-DeepONet with $\tanh(x)$, $\sin(x)$ and $\text{ReLU}(x)$ as activation functions obtains excellent predictions given limited training samples. However, the Causality-DeepONet with Sigmoid as activation functions is not convergent as expected. By shifting the Sigmoid
\begin{equation}
\sigma(x) = \dfrac{1}{1+e^{-x}} -\dfrac{1}{2}
    \label{eq: sigmoid-shift}
\end{equation}
could improve the results, as shown in Fig. \ref{Fig:Loss Causality DeepOnet-activation}(e) and case 5 in Table \ref{Table:Causality DeepONet-activation}.

\begin{table}[H]
\centering
\begin{threeparttable}
\caption{Relative $L_2$ Error for training and testing dataset when using Causality-DeepONet with different activation functions}
\label{Table:Causality DeepONet-activation}
\begin{tabular}{ccccc}
\hline
\multirow{2}{*}{Case} &  \multirow{2}{*}{Activation} & \multirow{2}{*}{Sample} &  \multicolumn{2}{c}{ Relative $L_2$ Error Eq. (\ref{Eq:RL2})}\\ \cline{4-5}
& & & Train & Test \\ \hline
1   & $\tanh(x)$           & 10 &0.012&0.015 \\
2   & $\sin(x)$            & 10 &0.013&0.019 \\
3   &  $\text{Sigmoid}(x)$  & 10 &0.165&0.143\\
4   &  $\text{ReLU}(x)$    & 10 &0.008&0.015 \\
5   & $\text{Custom Sigmoid  } (\text{Eq. } (\ref{eq: sigmoid-shift}))$ & 10 & 0.02&0.024\\
\hline
\end{tabular}
\end{threeparttable}
\end{table}

From the above discussion it could be observed that the proposed Causality-DeepONet is able to predict the response of the problem considered with a good accuracy. We also study the effect of training samples in the accuracy of predicted response. For this purpose we consider different samples with same network size and other hyper-parameters. The statistical properties of the different training samples are shown in Table \ref{table:data stats} in Appendix \ref{Appendix: data statistics}.
It could be noted that the training samples in datasets Train-I, Train-II, Train-III are exclusively different. The training samples in datasets Train-II and Train-III are included in the dataset Train-IV. The training samples in dataset Train-I and Train-IV are included in dataset Train-V. The training samples in dataset Train-V are included in dataset Train-VI.

\begin{figure}[H]
\centering
\subfigure{\includegraphics[width=0.32\textwidth]{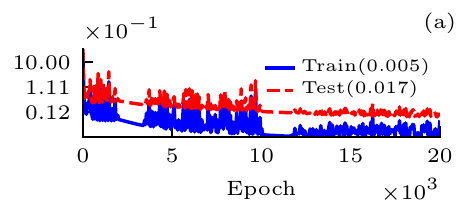}}
\subfigure{\includegraphics[width=0.32\textwidth]{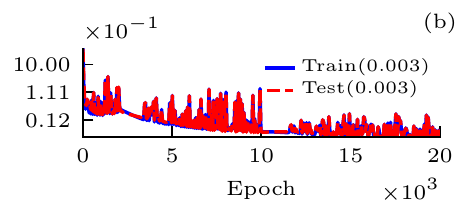}}
\subfigure{\includegraphics[width=0.32\textwidth]{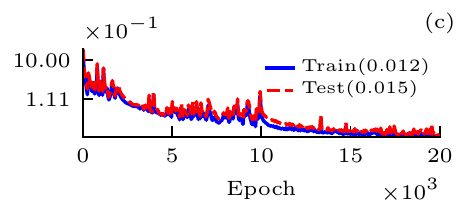}}
\subfigure{\includegraphics[width=0.32\textwidth]{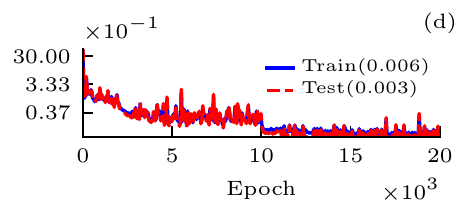}}
\subfigure{\includegraphics[width=0.32\textwidth]{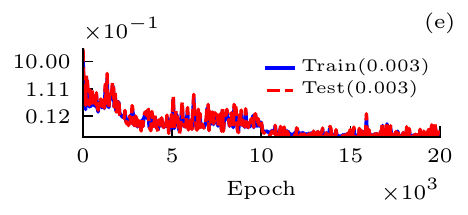}}
\subfigure{\includegraphics[width=0.32\textwidth]{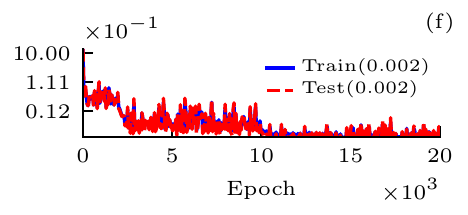}}
\caption{\textbf{Relative $\bm{L_2}$ Error for Training and Testing  Dataset when Using Causality-DeepONet with Different Number of Training Samples:} The plots (a)-(e) correspond to the cases Case-1 to Case-5 in Table \ref{Table:Causality DeepONet-nsample}, respectively.
}\label{Fig:Loss Causality DeepOnet-nsample}
\end{figure}

We consider a network size of [4000]-[120]$\times$2-[120] for branch and [1]-[120]$\times$2-[120] for trunk. The activation function considered is $\tanh(x)$. To avoid overfitting, we consider $L_2$ weight regularization with a coefficient $1\times 10^{-4}$ for the  branch net. The learning rate considered is $10^{-3}$ in the first 2000 epochs, then $10^{-4}$ in the 2000 to 10000 epochs, and $10^{-5}$ for the 10000 to 20000 epochs. The relative $L_2$ errors of different cases with different sample in training are shown in Table \ref{Table:Causality DeepONet-nsample} and Fig. \ref{Fig:Loss Causality DeepOnet-nsample}. It could be observed that the $L_2$ error is small even with smaller number of training set and the accuracy increases with the increase in number of training set, though the improvements in accuracy is limited with increase in number of samples.
On the other hand, it could also be observed that the performance of Causality-DeepONet of Case-3 in Table \ref{Table:Causality DeepONet-nsample} with 10 training samples that have larger deviation is poor comparing with Case-4 to Case-6 in Table \ref{Table:Causality DeepONet-nsample}. The further observation from Table \ref{Table:Causality DeepONet-nsample} is that the training relative $L_2$ error is greater than the testing relative $L_2$ error in Case-4 of Table \ref{Table:Causality DeepONet-nsample}, given the fact that dataset Train-IV contains Train-II and Train-III.

\begin{table}[H]
\centering
\begin{threeparttable}
\caption{Relative ${L_2}$ error for training and testing dataset when using Causality-DeepONet with different numbers of training samples}
\label{Table:Causality DeepONet-nsample}
\begin{tabular}{cccc}
\hline
\multirow{2}{*}{Case} &   \multirow{2}{*}{Sample} &  \multicolumn{2}{c}{Relative $L_2$ Error Eq. (\ref{Eq:RL2})} \\ \cline{3-4}
&  & Train & Test \\ \hline
1   &   7  &0.005&0.017 \\
2   &   8  &0.003&0.003 \\
3   &   10  &0.012&0.015 \\
4   &   20  &0.006 &0.003 \\
5   &   50 &0.003&0.003 \\
6   &   100 &0.002&0.002 \\
\hline
\end{tabular}
\end{threeparttable}
\end{table}

As discussed in section \ref{Section:Problem statement}, the proposed Causality-DeepONet involves both the phenomenon of causality and convolution. To evaluate the importance of convolution on the accuracy of prediction, we study the method only with causality but without convolution. The neural network considered for this purpose has form
\begin{equation}
    \mathcal{R}_{c}(\ddot{u}_{g})(t) \sim \sum_{k=1}^{N}\sum_{i=1}^{M}c_{i}^{k}\sigma
    _{b}\left(\sum_{j=1}^{\left\lfloor \frac{t}{h}\right\rfloor }W_{i,j}^{k}\ddot{u}_{g}\left(  s_{j}\right) +\sum_{j=\left\lfloor \frac{t}{h}\right\rfloor +1}^{\left\lfloor \frac{T}{h}\right\rfloor }W_{i,j}^{k}0 + B_{i}^{k}
\right)  \sigma_{t}\left(  \bm{w}_{k} t+b_{k}\right)
\label{eq: causality-deepONet-noconv}
\end{equation}

The difference between the formulation given by Eq. (\ref{eq: causality-deepONet-noconv}) and the proposed Causality-DeepONet (Eq. (\ref{eq: UAT_operator_causal})) is in the difference in the weights of the branch. In the case of the formulation without convolution, the weights are $W_{i,j}^{k}$, whereas in the case of the proposed Causality-DeepONet (with convolution) the weights are $ W_{i,\left\lfloor \frac{T}{h} \right\rfloor -\left\lfloor \frac{t}{h}\right\rfloor +j}^{k}$.

\begin{figure}[H]
\centering
\subfigure{\includegraphics[width=0.48\textwidth]{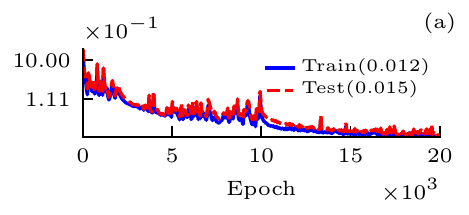}}
\subfigure{\includegraphics[width=0.48\textwidth]{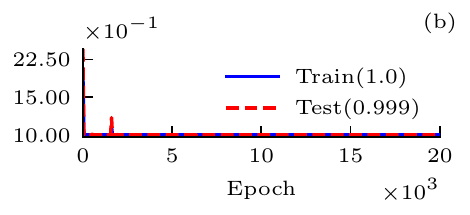}}

\caption{\textbf{Relative $\bm{L_2}$ Error for Training and Testing  Dataset when Using Causality-DeepONet with or without Convolutions:}
(a) Causality-DeepONet with convolution with loss function Eq. (\ref{Eq:loss deeponet 2}), 10 training samples.
(b) Proposed Net with causality only with loss function Eq. (\ref{Eq:loss deeponet 2}), 100 training samples.
}\label{Fig:Loss Causality DeepOnet-conv}
\end{figure}

 To study the effect of convolution, we compared the results of the problem with Causality-DeepONet and Causality-DeepONet without convolution.
A network size of $[4000]$-[120]$\times$2-$[120]$ for branch and  [1]-[120]$\times$2-[120] for the trunk are considered in both the cases. The activation function considered for both cases is $\tanh(x)$. To avoid overfitting, we consider $L_2$ weight regularization for the branch net with a coefficient $10^{-4}$. The learning rate considered is $10^{-3}$ in the first 2000 epochs, then $10^{-4}$ in the 2000 to 10000 epochs, and $10^{-5}$ for the rest of the epoch up to 20000 epochs. The case without convolution is considered to be trained with 100 training samples. However, the case with convolution is trained with 10 training samples only. Fig. \ref{Fig:Loss Causality DeepOnet-conv} shows the relative $L_2$ error for both the cases. It can be observed that the variant without convolution is not able to provide satisfactory accuracy.

\begin{figure}[H]
\def\svgwidth{18cm}
\centering  \import{Fig}{NonConvDeepONet.pdf_tex}
\caption{\textbf{Schematic Diagram of Causality-DeepONet without Convolution:} A schematic of the Causality-DeepONet without convolution showing branch and the trunk net along with the input data and output. The number of input neurons for the branch is equal to the number of sensor points in the input signals. The input signals for the branch, however, are replaced by a zero-padding for the future information.}
\label{fig:nonconv-cdeeponet}
\end{figure}
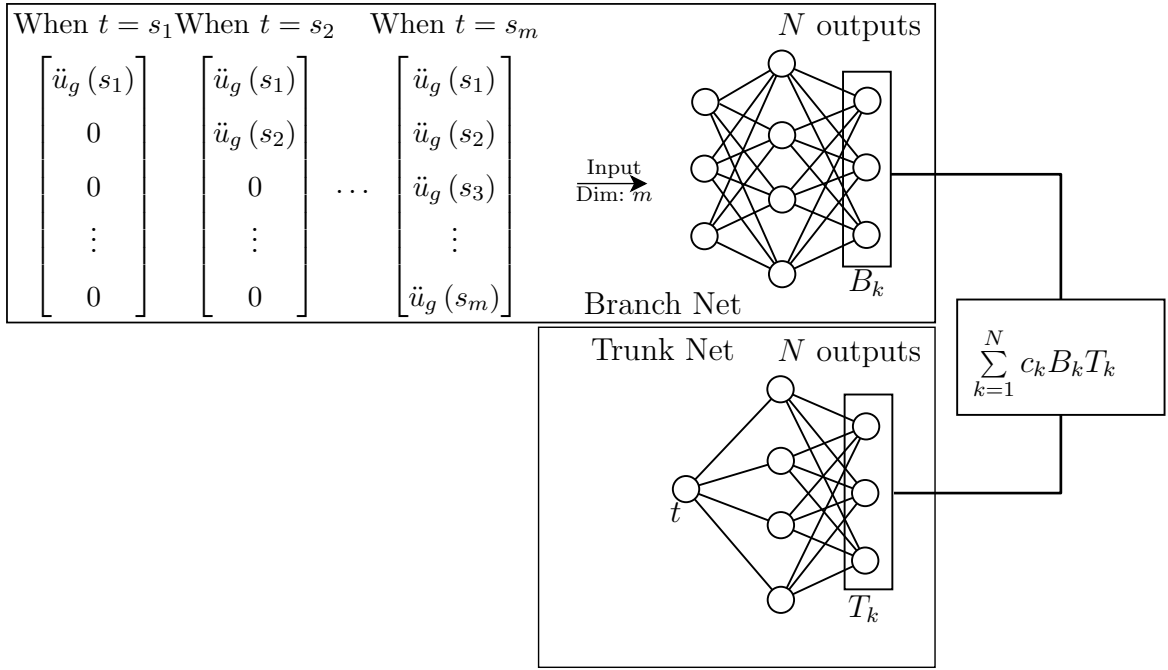

\section{Conclusion and future works}
\label{Section:Conclusions}
In this paper, we have studied how to improve the accuracy of DeepONet for causal oscillatory linear dynamical systems. Two new variants of DeepONet, the multi-scale DeepONet and the Causality-DeepONet are proposed. As an application,  we considered the problem of learning the mapping between earthquake ground accelerations and building's causal responses, which are both highly oscillatory. In the multi-scale DeepONet, multi-scale neural networks are used in the trunk net. Meanwhile, the Causality-DeepONet includes both causality and convolution as specific domain knowledge in its design. Though the multi-scale DeepONet improved the training of the seismic response operator, it failed to give satisfactory prediction results in the test cases. However, the Causality-DeepONet is able to provide accurate predictions in the test cases as well. We have also studied the effect of the size of networks, the number of training samples, and the type of activation functions on the accuracy of prediction of responses using Causality-DeepONet. It is found that the proposed Causality-DeepONet can provide good accuracy in the prediction of response of the problem considered. 

For future work, the Causality-DeepONet for nonlinear problems such as nonlinear dynamics, nonlinear electrical circuits etc, may be considered. Also future studies should include establishing a solid mathematical foundation for the Causality-DeepONet by extending the work of \cite{tianpingchen1995} to the proposed framework of the Causality-DeepONet.

\section*{Acknowledgement}
The authors like to thank Prof. George Em Karniadakis for bringing the attention of this research project to our attention and, Lu Lu as well,  for helpful discussions and assistance during this work. The work of W. Cai is supported by the US National Science Foundation grant DMS-2207449. The work of K. Nath is supported by OSD/AFOSR MURI grant FA9550-20-1-0358.

%\bibliographystyle{plain}
%\bibliography{Reference}

\newpage
\appendix
\textbf{Appendices}

\section{Additional Tables}
\label{Appendix: data statistics}
In Table \ref{table:data stats}, we show the statistical properties of the training and testing dataset. The total number of testing dataset considered is 44 earthquake ground accelerations. These test data are considered for all the numerical examples. It is also important to note that in the case of training, dataset Train-I, Train-II and Train-III are completely different dataset. Training dataset Train-I is not included in training dataset Train-II, similarly, training dataset Train-I and Train-II are not included in training dataset Train-III. However, training dataset Train-IV includes training dataset Train-III and training dataset Train-V includes training dataset Train-IV as well. The training dataset Train-VI includes all the training data I to V. The training dataset Train-VI is used for the training of DeepONet, DeepONet with Gaussian Normalization, POD-DeepONet, MSDeepONet. The training dataset Train-I to Train-VI are the data used for the discussion in section \ref{Subsection:Causality DeepONet}.
\begin{table}[H]
\centering
\caption{Properties of the earthquake records considered for training and testing of neural networks.}
\label{table:data stats}
\begin{tabular}{||l|c|c|c|c|c|c|c|c||}
\hline
 Case & &Test & Train-I &Train-II & Train-III& Train-IV & Train-V & Train-VI \\ \hline
 Samples&& 44 & 7 & 8 & 10 & 20 & 50 & 100 \\\hline
 \multirow{4}{*}{Mag.}
 & Min      &  4.30         &6.30       &4.90        & 4.92        &4.90       &4.90          & 4.70              \\
 &Max       &  7.90         & 7.62     & 7.62        & 7.62        & 7.62      & 7.62         & 7.62              \\
 &Mean      &  6.58         & 7.29     & 6.48        & 7.07        & 6.85      & 7.01         & 7.10              \\
 &SD        &  0.82         & 0.46     & 1.04        & 0.83        & 0.94      & 0.83         & 0.82              \\
 \hline
 \multirow{4}{*}{\begin{tabular}{@{}l@{}}
      PGA  \\
    $\max|\ddot{u}_g|$
 \end{tabular}
 }
 & Min      & 0.00103       & 0.00966  & 0.00499        & 0.00245     & 0.00245   & 0.00119      & 0.00119          \\
 &Max       & 0.35726       & 0.14425  & 0.11854        & 0.31313     & 0.31313   & 0.31313      & 0.31313          \\
 &Mean      & 0.05298       & 0.05915  & 0.05262        & 0.10973     & 0.08027   & 0.06972      & 0.07277          \\
 &SD        & 0.06428       & 0.04255  & 0.03863        & 0.11792     & 0.09269   & 0.07151      & 0.06338          \\
 \hline
 \multirow{4}{*}{Max $\ddot{u}_g$}
 & Min      & 0.00103       & 0.00966  & 0.00498        & 0.00245     & 0.00245   & 0.00118      & 0.00118          \\
 &Max       & 0.35726       & 0.14425  & 0.09805        & 0.30114     & 0.30114   & 0.30114      & 0.30114          \\
 &Mean      & 0.04883       & 0.05610  & 0.04460        & 0.10440     & 0.07439   & 0.06440      & 0.06661          \\
 &SD        & 0.06199       & 0.04148  & 0.03279        & 0.11204     & 0.08818   & 0.06795      & 0.06008          \\
 \hline
 \multirow{4}{*}{Min  $\ddot{u}_g$}
 & Min      & -0.27870     & -0.13713 & -0.11854          & -0.31313    & -0.31313  & -0.31313     & -0.31313             \\
 &Max       & -0.00095     & -0.00671 & -0.00499          & -0.00233    & -0.00233  & -0.00119     & -0.00112            \\
 &Mean      & -0.04820     & -0.05413 & -0.05262          & -0.09849    & -0.07365  & -0.06312     & -0.06791             \\
 &SD        &  0.05502     & 0.04177  & 0.03863        & 0.10559     & 0.08314   & 0.06391      &  0.05899           \\
 \hline
 \multirow{4}{*}{\begin{tabular}{@{}l@{}}
      Energy  \\
    $\|\ddot{u}_g\|^2$
 \end{tabular}}

 & Min      &  0.00020     & 0.00632  & 0.00188        & 0.00043     & 0.00043   & 0.00016      & 0.00016          \\
 &Max       &  3.22599     & 2.47209  & 2.05953        & 7.31542     & 7.31542   & 7.31542      & 7.31542          \\
 &Mean      &  0.38912     & 0.57996  & 0.44703        & 1.66238     & 1.03450   & 0.77436      & 0.82185          \\
 &SD        &  0.71029     & 0.81379  & 0.65613        & 2.40233     & 1.86043   & 1.34532      & 1.22943          \\
 \hline
\end{tabular}
\end{table}

\section{Additional Figures}
\subsection{A typical ground acceleration due to earthquake before and after processing}
\label{Appendix: processed signals}
\begin{figure}[H]
\centering
\subfigure{\includegraphics[width=0.8\textwidth]{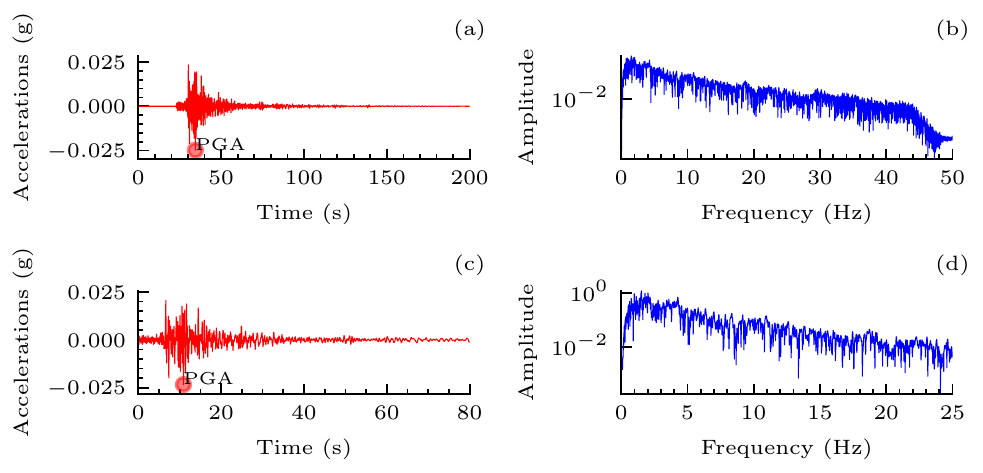}}
\caption{\textbf{Ground Acceleration:} The ground acceleration due to 14383980 earthquake recorded at station North Hollywood, 2008 (a) Time history of the acceleration. (b) Frequency spectrum. (c) The resampled acceleration. (d) The frequency spectrum of resampled acceleration.}
\label{Figure:earthquake-process}
\end{figure}

Fig. \ref{Figure:earthquake-process}(a)-(b) show ground acceleration due to earthquake 14383980 before and after processing. The ground acceleration record are recorded with $\delta t=0.005$ sec. The signal is passed through a butterworth filter. It is also important to note that the length of the signal is 200 second which is much higher than the other signal considered. Thus we have not considered initial 23.56 second acceleration and we have not considered the earthquake after 103.56 sec. The acceleration removed accounts for 0.7\% of total energy.
\subsection{Additional results for numerical study for DeepONet}
\label{Appendix: DeepONet-further-tests}
\begin{figure}[H]
\centering
\subfigure{\includegraphics[width=0.4\textwidth]{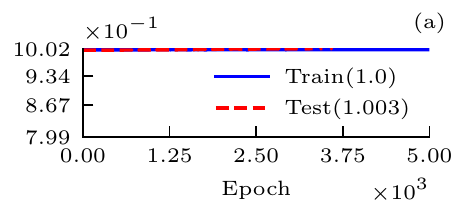}}
\subfigure{\includegraphics[width=0.4\textwidth]{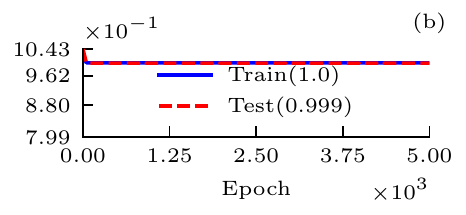}}
\subfigure{\includegraphics[width=0.4\textwidth]{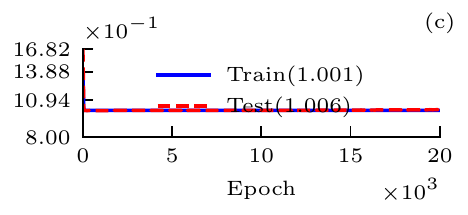}}
\subfigure{\includegraphics[width=0.4\textwidth]{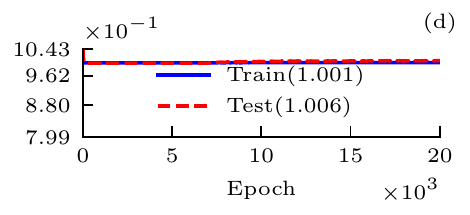}}
\caption{\textbf{Relative $\bm{L_2}$ Error for DeepONet:}  All of these cases are with loss function Eq. (\ref{Eq:loss deeponet 2}). The activation function for all the cases are $\text{ReLU}(x)$. To avoid overfitting, we consider using dropout rate 0.01 for the branch net and with $3\times 10^{-5}$ $L_2$ regularization. The network size is [4000]-[200]$\times$3-[200] and [1]-[200]$\times$3-[200] for the branch, the trunk net, respectively.
(a) In this case, $t$ is re-scaled to $[0,1]$, the learning rate is $10^{-4}$ in the first 1000 epochs, then $10^{-5}$ in the 1000 to 3000 epochs, at last $10^{-6}$ for the 3000 to 5000 epochs.
(b) In this case, the learning rate are set to be $10^{-4}$ for all the 5000 epochs.
(c) This is the case with same setting as case-3 in Table \ref{Table:DeepONet} but training up to 20000 epochs with same learning rate $10^{-4}$.
(d) This is the case with same setting as case-6 in Table \ref{Table:DeepONet} but training up to 20000 epochs with same learning rate $10^{-4}$.
}\label{Fig:Loss DeepOnet-extra}
\end{figure}
\begin{figure}[H]
\centering
\subfigure{\includegraphics[width=0.32\textwidth]{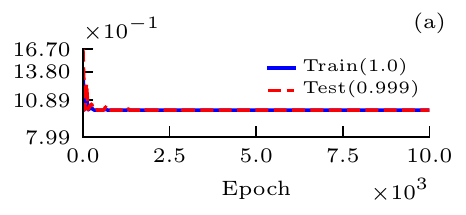}}
\subfigure{\includegraphics[width=0.32\textwidth]{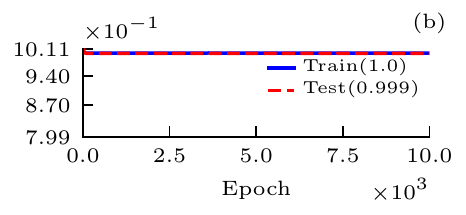}}
\subfigure{\includegraphics[width=0.32\textwidth]{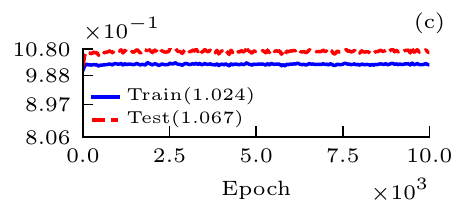}}
\caption{\textbf{Relative $\bm{L_2}$ Error for DeepONet with Different Activation Functions:}  All of these cases are with loss function Eq. (\ref{Eq:loss deeponet 2}) and are trained 10000 epochs. To avoid overfitting, we consider using dropout rate 0.01 for the branch net and with $3\times 10^{-5}$ $L_2$ regularization. The learning rate is $10^{-4}$ for all the cases. The network size is [4000]-[200]$\times$3-[200] and [1]-[200]$\times$3-[200] for the branch, the trunk net, respectively.
(a) with $\text{Sigmoid}(x)$ as activation function.
(b) with $\tanh(x)$ as activation function.
(c) with $\sin(x)$ as activation function.
}\label{Fig:Loss DeepOnet-extra-activation}
\end{figure}
\label{Appendix: predictions-DeepONet}
\begin{figure}[H]
\centering
\subfigure{\includegraphics[width=\textwidth]{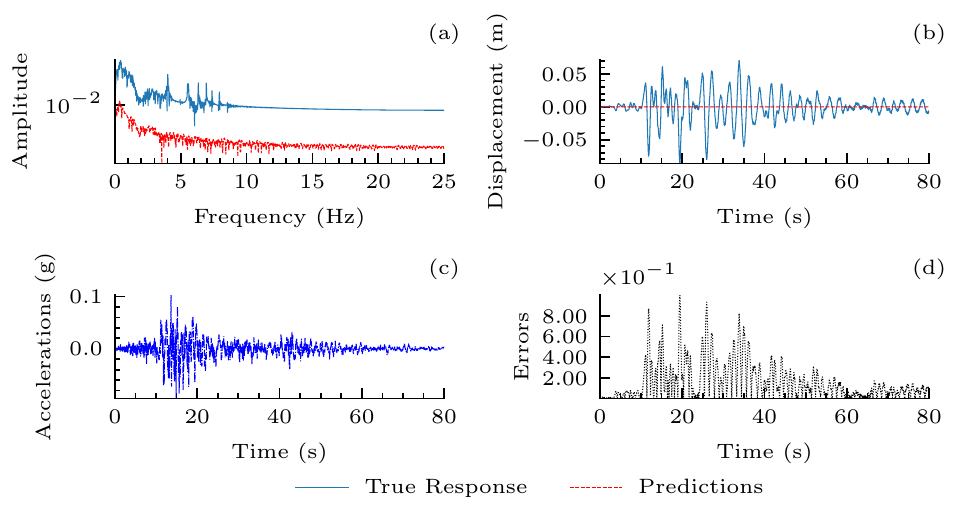}}
\caption{\textbf{One of the Training Samples with Predictions of DeepONet:} One of the training samples of DeepONet for Case-4 in Table. \ref{Table:DeepONet} with Loss Eq. (\ref{Eq:loss deeponet 2}).
(a) The Amplitude of the prediction and true response in Fourier Domain,
(b) The prediction and true response,
(c) The corresponding input signals,
(d) The relative error Eq. (\ref{Eq:RE}).
}
\label{Fig:DeepOnet Pred - Training}
\end{figure}

\newpage
\begin{figure}[H]
\centering
\subfigure{\includegraphics[width=\textwidth]{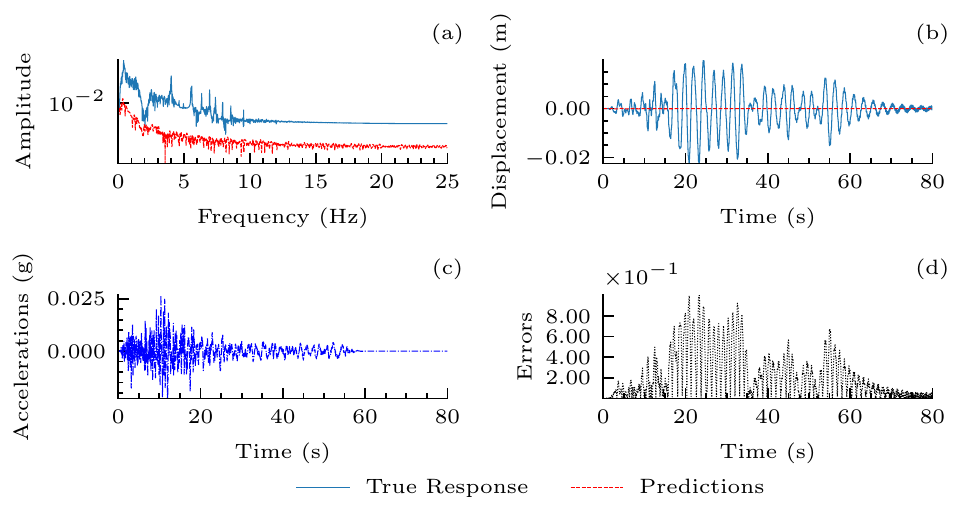}}
\caption{\textbf{The Best Case of the Predictions of DeepONet:} The best predictions for the DeepONet for testing cases for Case-4 in Table. \ref{Table:DeepONet} with Loss Eq. (\ref{Eq:loss deeponet 2}).
(a) The Amplitude of the prediction and true response in Fourier Domain,
(b) The prediction and true response,
(c) The corresponding input signals,
(d) The relative error Eq. (\ref{Eq:RE}).
}
\label{Fig:DeepOnet Pred-Best}
\end{figure}

\subsection{Additional results for numerical study for POD-DeepONet}
\label{Appendix: POD-DeepONet-further-tests}
\begin{figure}[h]
\centering
\subfigure{\includegraphics[width=0.32\textwidth]{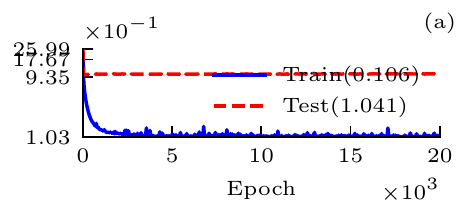}}
\subfigure{\includegraphics[width=0.32\textwidth]{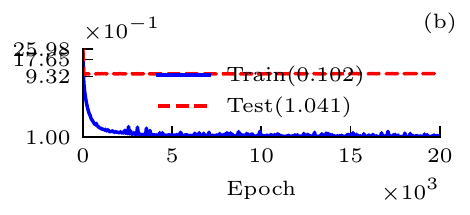}}
\subfigure{\includegraphics[width=0.32\textwidth]{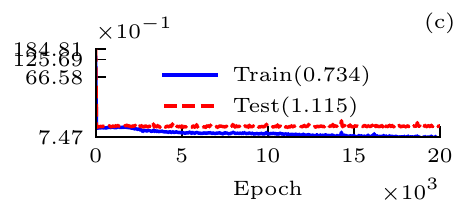}}
\caption{\textbf{Relative $\bm{L_2}$ Error for POD-DeepONet with Different Activation Functions:}  All of these cases are with loss function Eq. (\ref{Eq:loss deeponet 2}) and are trained 20000 epochs. To avoid overfitting, we consider using $10^{-6}$ $L_2$ regularization for the net. The learning rate is $10^{-4}$ for all the cases. The network size is [4000]-[200]$\times$3-[100].
(a) with $\sin(x)$ as activation function.
(b) with $\tanh(x)$ as activation function.
(c) with $\text{Sigmoid}(x)$ as activation function.
}\label{Fig:Loss POD-DeepOnet-extra}
\end{figure}

\newpage
\begin{figure}[H]
\centering
\subfigure{\includegraphics[width=\textwidth]{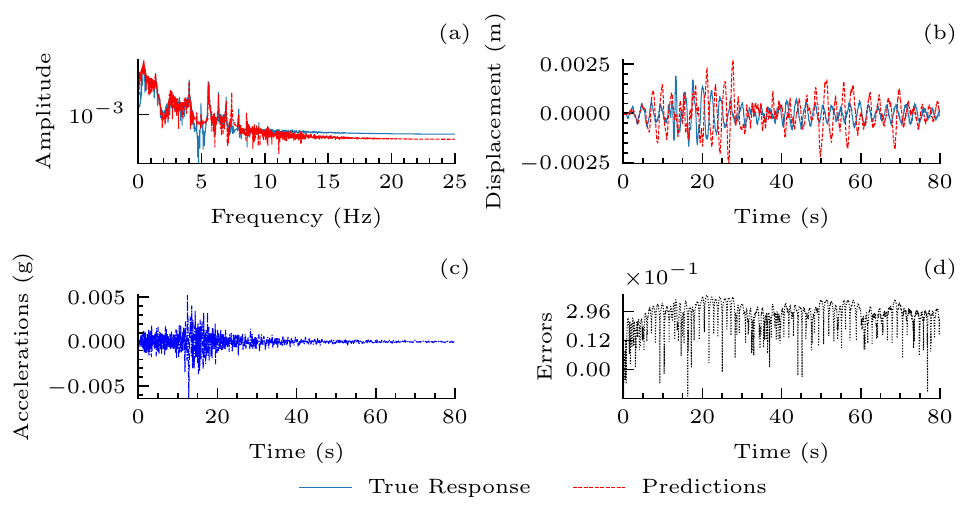}}

\caption{\textbf{One of the Training Samples with Predictions of POD-DeepONet:} One of the training samples of POD-DeepONet with Loss Eq. (\ref{Eq:loss deeponet 2}).
(a) The Amplitude of the prediction and true response in Fourier Domain,
(b) The prediction and true response,
(c) The corresponding input signals,
(d) The relative error Eq. (\ref{Eq:RE}).
}
\label{Fig:POD-DeepOnet Pred - Training}
\end{figure}

\begin{figure}[H]
\centering
\subfigure{\includegraphics[width=\textwidth]{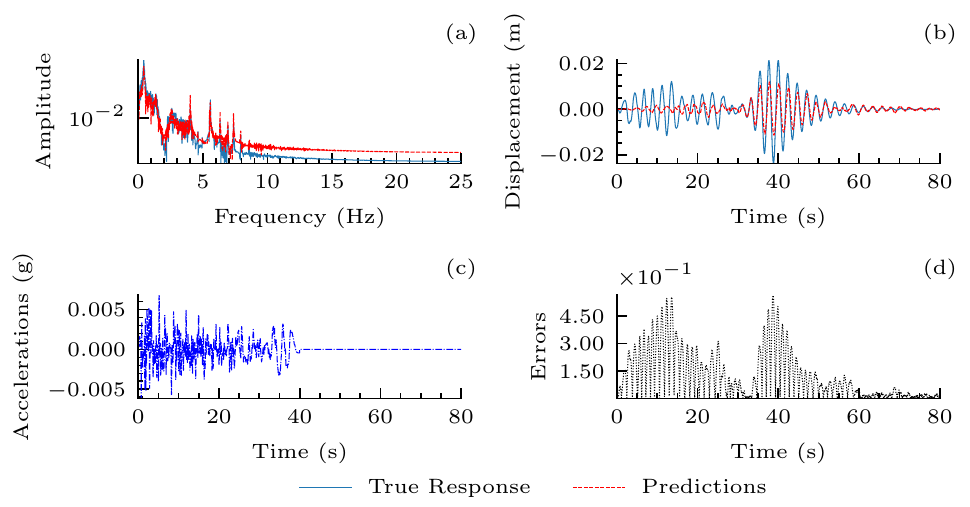}}

\caption{\textbf{The Best Case of the Predictions of POD-DeepONet:} The best predictions for the POD-DeepONet for testing cases with Loss Eq. (\ref{Eq:loss deeponet 2}).
(a) The Amplitude of the prediction and true response in Fourier Domain,
(b) The prediction and true response,
(c) The corresponding input signals,
(d) The relative error Eq. (\ref{Eq:RE}).
}
\label{Fig:POD-DeepOnet Pred-Best}
\end{figure}
\newpage
\subsection{Additional results for numerical study for Multi-scale DeepONet}
\begin{figure}[H]
\centering
\subfigure{\includegraphics[width=\textwidth]{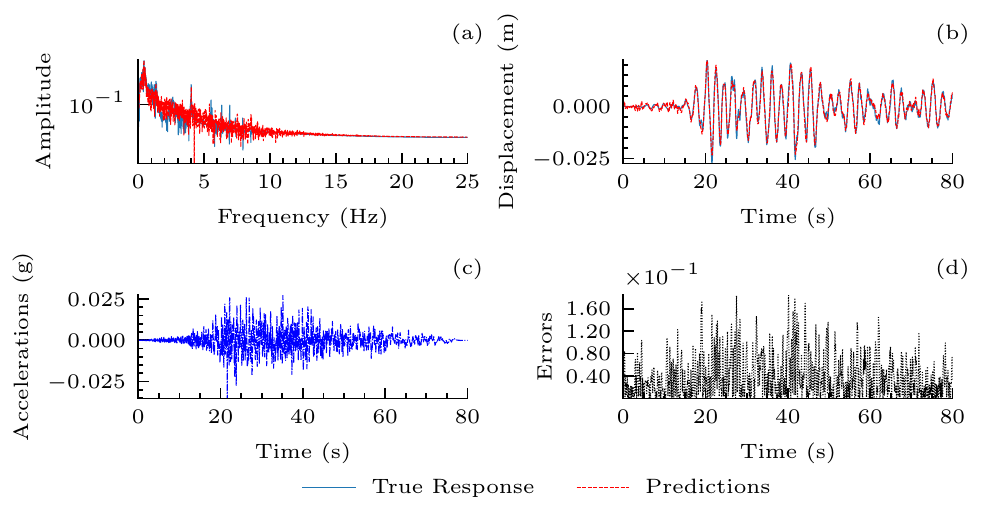}}

\caption{\textbf{One of the Training Samples with Predictions of MS-DeepONet:} One of the training samples of MS-DeepONet with Loss Eq. (\ref{Eq:loss deeponet 2}).
(a) The Amplitude of the prediction and true response in Fourier Domain,
(b) The prediction and true response,
(c) The corresponding input signals,
(d) The relative error Eq. (\ref{Eq:RE}).
}
\label{Fig:MSDeepOnet Pred - Training}
\end{figure}

\begin{figure}[H]
\centering
\subfigure{\includegraphics[width=\textwidth]{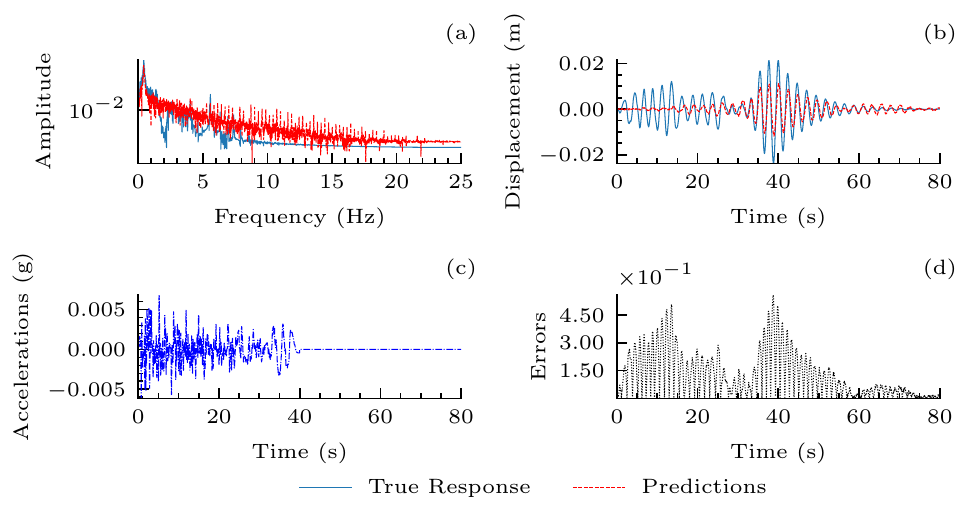}}

\caption{\textbf{The Best Case of the Predictions of MS-DeepONet:} The best predictions for the MS-DeepONet for testing cases with Loss Eq. (\ref{Eq:loss deeponet 2}).
(a) The Amplitude of the prediction and true response in Fourier Domain,
(b) The prediction and true response,
(c) The corresponding input signals,
(d) The relative error Eq. (\ref{Eq:RE}).
}
\label{Fig:MSDeepOnet Pred-Best}
\end{figure}

\end{document}

%% file: Fig/DeepONet.pdf_tex
%% Creator: Inkscape 1.2.1 (9c6d41e410, 2022-07-14), www.inkscape.org
%% Accompanies image file 'DeepOnet.pdf' (pdf, eps, ps)
\begingroup%
  \makeatletter%
  \providecommand\color[2][]{%
    \errmessage{(Inkscape) Color is used for the text in Inkscape, but the package 'color.sty' is not loaded}%
    \renewcommand\color[2][]{}%
  }%
  \providecommand\transparent[1]{%
    \errmessage{(Inkscape) Transparency is used (non-zero) for the text in Inkscape, but the package 'transparent.sty' is not loaded}%
    \renewcommand\transparent[1]{}%
  }%
  \providecommand\rotatebox[2]{#2}%
  \newcommand*\fsize{\dimexpr\f@size pt\relax}%
  \newcommand*\lineheight[1]{\fontsize{\fsize}{#1\fsize}\selectfont}%
  \ifx\svgwidth\undefined%
    \setlength{\unitlength}{533.25bp}%
    \ifx\svgscale\undefined%
      \relax%
    \else%
      \setlength{\unitlength}{\unitlength * \real{\svgscale}}%
    \fi%
  \else%
    \setlength{\unitlength}{\svgwidth}%
  \fi%
  \global\let\svgwidth\undefined%
  \global\let\svgscale\undefined%
  \makeatother%
  \begin{picture}(1,0.52039381)%
    \lineheight{1}%
    \setlength\tabcolsep{0pt}%
    \put(0,0){\includegraphics[width=\unitlength,page=1]{DeepOnet.pdf}}%
    \put(0.62374394,0.24808502){\color[rgb]{0,0,0}\makebox(0,0)[t]{\lineheight{1.25}\smash{\begin{tabular}[t]{c}$N$ outputs \end{tabular}}}}%
    \put(0.62374394,0.48803616){\color[rgb]{0,0,0}\makebox(0,0)[t]{\lineheight{1.25}\smash{\begin{tabular}[t]{c}$N$ outputs \end{tabular}}}}%
    \put(0,0){\includegraphics[width=\unitlength,page=2]{DeepOnet.pdf}}%
    \put(0.45186053,0.38557558){\color[rgb]{0,0,0}\makebox(0,0)[t]{\lineheight{1.25}\smash{\begin{tabular}[t]{c}\scriptsize Input \end{tabular}}}}%
    \put(0.45186053,0.36557558){\color[rgb]{0,0,0}\makebox(0,0)[t]{\lineheight{1.25}\smash{\begin{tabular}[t]{c}\scriptsize Dim: $m$\end{tabular}}}}%
    \put(0,0){\includegraphics[width=\unitlength,page=3]{DeepOnet.pdf}}%
    \put(0.48752677,0.2821552){\color[rgb]{0,0,0}\makebox(0,0)[t]{\lineheight{1.25}\smash{\begin{tabular}[t]{c}Branch Net\end{tabular}}}}%
    \put(0,0){\includegraphics[width=\unitlength,page=4]{DeepOnet.pdf}}%
    \put(0.48752677,0.2483481){\color[rgb]{0,0,0}\makebox(0,0)[t]{\lineheight{1.25}\smash{\begin{tabular}[t]{c}Trunk Net\end{tabular}}}}%
    \put(0,0){\includegraphics[width=\unitlength,page=5]{DeepOnet.pdf}}%
    \put(0.7143896,0.23985795){\color[rgb]{0,0,0}\makebox(0,0)[lt]{\lineheight{4.25}\smash{$\sum\limits_{k=1}^{N}c_{k}{B_k}{T_k}$}}}%
    \put(0.00482819,0.39121552){\color[rgb]{0,0,0}\makebox(0,0)[lt]{\lineheight{1.25}\smash{
    {\footnotesize
    \begin{tabular}{cccc}   
    When $x = x_1$ & When $x = x_2$ &                   & When $x = x_m$ \\
                   $\begin{bmatrix} 
                   f\left( x_1 \right)\\
                   f\left( x_2 \right)\\
                    \vdots \\
                    f\left( x_{m-1} \right)\\  
                    f\left( x_{m} \right)
                   \end{bmatrix}$ 
                   &
                    $\begin{bmatrix} 
                    f\left( x_1 \right)\\
                    f\left( x_2 \right)\\
                    \vdots \\
                    f\left( x_{m-1} \right)\\  
                    f\left( x_{m} \right)
                    \end{bmatrix}$ 
                   & $\ldots$ &
                    $\begin{bmatrix} 
                    f\left( x_1 \right)\\
                    f\left( x_2 \right)\\
                    \vdots \\
                    f\left( x_{m-1} \right)\\  
                    f\left( x_{m} \right)
                    \end{bmatrix}$ 
    \end{tabular}}
    }}}%
    \put(0,0){\includegraphics[width=\unitlength,page=6]{DeepOnet.pdf}}%
    \put(0.49280564,0.12901714){\color[rgb]{0,0,0}\makebox(0,0)[lt]{\lineheight{1.25}\smash{\begin{tabular}[t]{l}$x $\end{tabular}}}}%
    \put(0.62280564,0.29901714){\color[rgb]{0,0,0}\makebox(0,0)[lt]{\lineheight{1.25}\smash{\begin{tabular}[t]{l}$B_k$\end{tabular}}}}%
    \put(0.62280564,0.05901714){\color[rgb]{0,0,0}\makebox(0,0)[lt]{\lineheight{1.25}\smash{\begin{tabular}[t]{l}$T_k$\end{tabular}}}}%
  \end{picture}%
\endgroup%

%% file: Fig/CausalityDeepONet.pdf_tex
%% Creator: Inkscape 1.2.1 (9c6d41e410, 2022-07-14), www.inkscape.org
%% Accompanies image file 'DeepOnet.pdf' (pdf, eps, ps)
\begingroup%
  \makeatletter%
  \providecommand\color[2][]{%
    \errmessage{(Inkscape) Color is used for the text in Inkscape, but the package 'color.sty' is not loaded}%
    \renewcommand\color[2][]{}%
  }%
  \providecommand\transparent[1]{%
    \errmessage{(Inkscape) Transparency is used (non-zero) for the text in Inkscape, but the package 'transparent.sty' is not loaded}%
    \renewcommand\transparent[1]{}%
  }%
  \providecommand\rotatebox[2]{#2}%
  \newcommand*\fsize{\dimexpr\f@size pt\relax}%
  \newcommand*\lineheight[1]{\fontsize{\fsize}{#1\fsize}\selectfont}%
  \ifx\svgwidth\undefined%
    \setlength{\unitlength}{533.25bp}%
    \ifx\svgscale\undefined%
      \relax%
    \else%
      \setlength{\unitlength}{\unitlength * \real{\svgscale}}%
    \fi%
  \else%
    \setlength{\unitlength}{\svgwidth}%
  \fi%
  \global\let\svgwidth\undefined%
  \global\let\svgscale\undefined%
  \makeatother%
  \begin{picture}(1,0.52039381)%
    \lineheight{1}%
    \setlength\tabcolsep{0pt}%
    \put(0,0){\includegraphics[width=\unitlength,page=1]{DeepOnet.pdf}}%
    \put(0.62374394,0.24808502){\color[rgb]{0,0,0}\makebox(0,0)[t]{\lineheight{1.25}\smash{\begin{tabular}[t]{c}$N$ outputs \end{tabular}}}}%
    \put(0.62374394,0.48803616){\color[rgb]{0,0,0}\makebox(0,0)[t]{\lineheight{1.25}\smash{\begin{tabular}[t]{c}$N$ outputs \end{tabular}}}}%
    \put(0,0){\includegraphics[width=\unitlength,page=2]{DeepOnet.pdf}}%
    \put(0.45186053,0.38557558){\color[rgb]{0,0,0}\makebox(0,0)[t]{\lineheight{1.25}\smash{\begin{tabular}[t]{c}\scriptsize Input \end{tabular}}}}%
    \put(0.45186053,0.36557558){\color[rgb]{0,0,0}\makebox(0,0)[t]{\lineheight{1.25}\smash{\begin{tabular}[t]{c}\scriptsize Dim: $m$\end{tabular}}}}%
    \put(0,0){\includegraphics[width=\unitlength,page=3]{DeepOnet.pdf}}%
    \put(0.48752677,0.2821552){\color[rgb]{0,0,0}\makebox(0,0)[t]{\lineheight{1.25}\smash{\begin{tabular}[t]{c}Branch Net\end{tabular}}}}%
    \put(0,0){\includegraphics[width=\unitlength,page=4]{DeepOnet.pdf}}%
    \put(0.48752677,0.2483481){\color[rgb]{0,0,0}\makebox(0,0)[t]{\lineheight{1.25}\smash{\begin{tabular}[t]{c}Trunk Net\end{tabular}}}}%
    \put(0,0){\includegraphics[width=\unitlength,page=5]{DeepOnet.pdf}}%
    \put(0.7143896,0.23985795){\color[rgb]{0,0,0}\makebox(0,0)[lt]{\lineheight{4.25}\smash{$\sum\limits_{k=1}^{N}c_{k}{B_k}{T_k}$}}}%
    \put(0.00482819,0.39121552){\color[rgb]{0,0,0}\makebox(0,0)[lt]{\lineheight{1.25}\smash{
    {\footnotesize
    \begin{tabular}{cccc}
    When $t = s_1$ & When $t = s_2$ &                   & When $t = s_m$ \\
    $ \begin{bmatrix} 
                    0\\
                    \vdots \\
                    0\\
                    0\\  
                    \ddot{u}_{g}\left( s_{1} \right)
                   \end{bmatrix} 
                     $
                   &
     $ \begin{bmatrix} 
                    0\\
                    \vdots \\
                    0\\
                    \ddot{u}_{g}\left( s_{1} \right)\\  
                    \ddot{u}_{g}\left( s_{2} \right)
                   \end{bmatrix} 
                     $
                   &$\ldots $&
                   $ \begin{bmatrix} 
                    \ddot{u}_{g}\left( s_1 \right)\\
                    \vdots \\
                    \ddot{u}_{g}\left( s_{m-2} \right)\\
                    \ddot{u}_{g}\left( s_{m-1} \right)\\  
                    \ddot{u}_{g}\left( s_{m} \right)
                   \end{bmatrix} 
                     $
    \end{tabular}}
    }}}%
    \put(0,0){\includegraphics[width=\unitlength,page=6]{DeepOnet.pdf}}%
    \put(0.49280564,0.12901714){\color[rgb]{0,0,0}\makebox(0,0)[lt]{\lineheight{1.25}\smash{\begin{tabular}[t]{l}$t $\end{tabular}}}}%
      \put(0.62280564,0.29901714){\color[rgb]{0,0,0}\makebox(0,0)[lt]{\lineheight{1.25}\smash{\begin{tabular}[t]{l}$B_k$\end{tabular}}}}%
    \put(0.62280564,0.05901714){\color[rgb]{0,0,0}\makebox(0,0)[lt]{\lineheight{1.25}\smash{\begin{tabular}[t]{l}$T_k$\end{tabular}}}}%
  \end{picture}%
\endgroup%

%% file: Fig/NonConvDeepONet.pdf_tex
%% Creator: Inkscape 1.2.1 (9c6d41e410, 2022-07-14), www.inkscape.org
%% Accompanies image file 'DeepOnet.pdf' (pdf, eps, ps)
\begingroup%
  \makeatletter%
  \providecommand\color[2][]{%
    \errmessage{(Inkscape) Color is used for the text in Inkscape, but the package 'color.sty' is not loaded}%
    \renewcommand\color[2][]{}%
  }%
  \providecommand\transparent[1]{%
    \errmessage{(Inkscape) Transparency is used (non-zero) for the text in Inkscape, but the package 'transparent.sty' is not loaded}%
    \renewcommand\transparent[1]{}%
  }%
  \providecommand\rotatebox[2]{#2}%
  \newcommand*\fsize{\dimexpr\f@size pt\relax}%
  \newcommand*\lineheight[1]{\fontsize{\fsize}{#1\fsize}\selectfont}%
  \ifx\svgwidth\undefined%
    \setlength{\unitlength}{533.25bp}%
    \ifx\svgscale\undefined%
      \relax%
    \else%
      \setlength{\unitlength}{\unitlength * \real{\svgscale}}%
    \fi%
  \else%
    \setlength{\unitlength}{\svgwidth}%
  \fi%
  \global\let\svgwidth\undefined%
  \global\let\svgscale\undefined%
  \makeatother%
  \begin{picture}(1,0.52039381)%
    \lineheight{1}%
    \setlength\tabcolsep{0pt}%
    \put(0,0){\includegraphics[width=\unitlength,page=1]{DeepOnet.pdf}}%
    \put(0.62374394,0.24808502){\color[rgb]{0,0,0}\makebox(0,0)[t]{\lineheight{1.25}\smash{\begin{tabular}[t]{c}$N$ outputs \end{tabular}}}}%
    \put(0.62374394,0.48803616){\color[rgb]{0,0,0}\makebox(0,0)[t]{\lineheight{1.25}\smash{\begin{tabular}[t]{c}$N$ outputs \end{tabular}}}}%
    \put(0,0){\includegraphics[width=\unitlength,page=2]{DeepOnet.pdf}}%
    \put(0.45186053,0.38557558){\color[rgb]{0,0,0}\makebox(0,0)[t]{\lineheight{1.25}\smash{\begin{tabular}[t]{c}\scriptsize Input \end{tabular}}}}%
    \put(0.45186053,0.36557558){\color[rgb]{0,0,0}\makebox(0,0)[t]{\lineheight{1.25}\smash{\begin{tabular}[t]{c}\scriptsize Dim: $m$\end{tabular}}}}%
    \put(0,0){\includegraphics[width=\unitlength,page=3]{DeepOnet.pdf}}%
    \put(0.48752677,0.2821552){\color[rgb]{0,0,0}\makebox(0,0)[t]{\lineheight{1.25}\smash{\begin{tabular}[t]{c}Branch Net\end{tabular}}}}%
    \put(0,0){\includegraphics[width=\unitlength,page=4]{DeepOnet.pdf}}%
    \put(0.48752677,0.2483481){\color[rgb]{0,0,0}\makebox(0,0)[t]{\lineheight{1.25}\smash{\begin{tabular}[t]{c}Trunk Net\end{tabular}}}}%
    \put(0,0){\includegraphics[width=\unitlength,page=5]{DeepOnet.pdf}}%
    \put(0.7143896,0.23985795){\color[rgb]{0,0,0}\makebox(0,0)[lt]{\lineheight{4.25}\smash{$\sum\limits_{k=1}^{N}c_{k}{B_k}{T_k}$}}}%
    \put(0.00482819,0.39121552){\color[rgb]{0,0,0}\makebox(0,0)[lt]{\lineheight{1.25}\smash{
    {\footnotesize
    \begin{tabular}{cccc}
    When $t = s_1$ & When $t = s_2$ &                   & When $t = s_m$ \\
    $ \begin{bmatrix} 
                    \ddot{u}_{g}\left( s_1 \right)\\
                    0\\
                    0\\
                    \vdots \\
                    0 
                   \end{bmatrix} 
                     $
                   &
     $ \begin{bmatrix} 
                    \ddot{u}_{g}\left( s_1 \right)\\
                    \ddot{u}_{g}\left( s_2 \right)\\
                    0\\
                    \vdots \\
                    0
                   \end{bmatrix} 
                     $
                   &$\ldots $&
                   $ \begin{bmatrix} 
                    \ddot{u}_{g}\left( s_1 \right)\\
                    \ddot{u}_{g}\left( s_2 \right)\\
                    \ddot{u}_{g}\left( s_{3} \right)\\  
                    \vdots \\
                    \ddot{u}_{g}\left( s_{m} \right)
                   \end{bmatrix} 
                     $
    \end{tabular}}
    }}}%
    \put(0,0){\includegraphics[width=\unitlength,page=6]{DeepOnet.pdf}}%
    \put(0.49280564,0.12901714){\color[rgb]{0,0,0}\makebox(0,0)[lt]{\lineheight{1.25}\smash{\begin{tabular}[t]{l}$t $\end{tabular}}}}%
       \put(0.62280564,0.29901714){\color[rgb]{0,0,0}\makebox(0,0)[lt]{\lineheight{1.25}\smash{\begin{tabular}[t]{l}$B_k$\end{tabular}}}}%
    \put(0.62280564,0.05901714){\color[rgb]{0,0,0}\makebox(0,0)[lt]{\lineheight{1.25}\smash{\begin{tabular}[t]{l}$T_k$\end{tabular}}}}% 
\end{picture}%
\endgroup%